\documentclass[12pt]{article} 

\usepackage{standalone}
\usepackage[margin=1in]{geometry}
\usepackage[utf8]{inputenc}
\usepackage[dvipsnames]{xcolor}
\usepackage[shortlabels]{enumitem}
\usepackage{graphicx}
\usepackage{amsmath,amssymb,amsthm,mathtools}
\usepackage{MnSymbol}
\usepackage{latexsym, verbatim}
\usepackage{appendix}
\usepackage{units}
\usepackage{amsmath}
\usepackage[normalem]{ulem}
\usepackage[hidelinks]{hyperref}
\usepackage{hyperref} 
\usepackage{soul}
\usepackage{nicematrix}
\usepackage{pgfplots}
\usepackage[ruled,linesnumbered]{algorithm2e}
\pgfplotsset{compat=1.18}
\usetikzlibrary{pgfplots.groupplots}
\usepackage{tikz-cd}

\mathtoolsset{showonlyrefs}
\newtheorem{theorem}{Theorem}[section]
\newtheorem{corollary}[theorem]{Corollary} 
\newtheorem{lemma}[theorem]{Lemma}
\newtheorem{proposition}[theorem]{Proposition}

\theoremstyle{definition}
\newtheorem{example}[theorem]{Example}
\newtheorem{definition}[theorem]{Definition}
\newtheorem{remark}[theorem]{Remark}

\newcommand{\Acal}{\mathcal{A}}
\newcommand{\Ucal}{\mathcal{U}}

\newcommand{\Hcal}{\mathcal{H}}

\newcommand{\Ccal}{\mathcal{C}}

\newcommand{\Vcal}{\mathcal{V}}

\newcommand{\CC}{\mathbb{C}}

\newcommand{\PP}{\mathbb{P}}
\newcommand{\RR}{\mathbb{R}}

\DeclareMathOperator{\Conv}{Conv}

\DeclareMathOperator{\codim}{codim}
\DeclareMathOperator{\Gr}{Gr}
\DeclareMathOperator{\CH}{\rm CH}

\DeclareMathOperator{\Sing}{Sing}

\newcommand{\bP}{\mathbb{P}}
\newcommand{\cA}{\mathcal{A}}

\title{{\bf Avoidance Loci of Real Projective Varieties}}
\author{Elizabeth Pratt\footnote{Department of Mathematics, 
University of California, Berkeley \texttt{epratt@berkeley.edu}} \ and Kexin Wang\footnote{Harvard John A. Paulson School of Engineering and Applied Sciences \texttt{kexin\_wang@g.harvard.edu}}}
\date{}

\begin{document}
\maketitle{}

\begin{abstract}
We study real linear spaces in projective space that avoid the real points of a non-degenerate projective variety. For a variety $X \subset \mathbb{P}^{n-1}$ with a real smooth point, we define the \emph{avoidance locus} $\mathcal{A}_k(X)$ as the subset of the real Grassmannian $\mathrm{Gr}(k,n)_{\mathbb{R}}$ consisting of linear spaces that meet $X$ transversely but contain no real point of $X$. Our construction generalizes the cone of positive polynomials on $\RR^n.$ We prove that the avoidance locus is an open semi-algebraic set equal to a union of regions in the complement of a higher Chow form, and that distinct regions are non-adjacent. We present explicit examples for linear spaces, curves, and surfaces, and provide bounds on the number of connected components of $\mathcal{A}_{n-1}(X)$ in terms of the topology of the real locus $X_{\mathbb{R}}$. Finally, we prove that avoidance loci are slice-convex.
\end{abstract}

\section{Introduction}

In this paper, we study real linear spaces that avoid the real points of a real projective variety. 
Throughout the paper, we consider a smooth $d$-dimensional complex variety $X\subseteq \PP^{n-1}$ which is non-degenerate (i.e. not contained in any hyperplane) and \emph{totally real}.

\begin{definition}[Totally real variety]
A complex variety $X$ is \emph{totally real} if it is defined over the real numbers and its real locus is Zariski dense.
\end{definition}

This terminology agrees with \cite{blekherman}. If the variety $X$ is defined over $\mathbb{R}$ and contains a real smooth point, then it is automatically totally real \cite[Proposition 7.6.2]{bochnak2013real}.
We denote the real locus of $X$ by $X_\RR$ and the Grassmannian of $(k-1)$-dimensional linear spaces in $\PP^{n-1}$ by $\Gr(k,n)$. By dimension we always mean projective dimension.

General linear spaces of projective dimension at most $n-d-2$ do not meet $X$. 
General linear spaces of dimension $n-d-1$ meet $X$ in finitely many points. 
If $\deg X$ is even, it could happen that all these points are complex. 
Generic linear spaces of dimension larger than $n-d-1$ meet $X$ in a positive-dimensional subvariety, but it could also happen that this subvariety contains no real point. 
We study the locus in the Grassmannian of the linear spaces where this occurs. 

\begin{definition}
Let $X\subseteq \PP^{n-1}$ be a smooth non-degenerate totally real variety of dimension $d$. The k-th avoidance locus is 
$$
\Acal_k(X) := \{ V \in \Gr(k,n)_\RR \, : \, V \cap X \text{ is transverse and } V\cap X_\RR = \varnothing\}.
$$
\end{definition}

\begin{remark}
    We emphasize that the avoidance locus consists of linear spaces which avoid the variety \emph{transversely}. For example, lines which meet a space curve in complex points are \emph{not} in the avoidance locus. Similarly, lines which avoid the real part of a plane curve but tangent to it at a pair of complex conjugate points are not in the avoidance locus. This differs from the convention of \cite{kaihnsa2019sixty}, who define the avoidance locus of a plane curve as the set of lines which do not intersect it.
\end{remark}

For example, when $k=1$ we have that $\Acal_k(X) = \PP^{n-1}_\RR \setminus X_\RR.$ When $k = n-d$, the set $\Acal_k(X)$ contains linear spaces meeting $X$ in finitely many complex points but no real points. 

Avoidance loci specialize to known constructions. 
When $X$ is a finite number of points, $\Acal_{n-1}(X)$ consists of all hyperplanes which do not pass through any of the points. More generally for a projective variety $X,$ the avoidance locus $\Acal_{n-1}(X)$  of hyperplanes is a union of open regions in the complement of the dual variety $X^\vee$. It is non-empty if and only if the variety is fully contained in some affine chart. Such varieties are called \emph{very compact}. 

\begin{example}[The Trott curve]
    The Trott curve is a real quartic curve with four real non-nested components and $28$ bitangents, all of which are real. Its equation is
    \[144(x^4+ y^4) - 225z^2(x^2+y^2) + 350x^2y^2 + 81z^4 = 0.\]
    Its dual curve has $28$ monomials of degree $12$. Its complement in $(\bP^2)^\vee$ has $29$ components, which we computed using \texttt{HypersurfaceRegions.jl}~\cite{breiding2025computing}. Seven of these make up the avoidance locus. The right side of Figure \ref{fig:trott} shows $25$ of these regions in the affine chart given by $z = 1$, with the avoidance locus shaded. Note that the four unbounded regions are actually two regions in projective space.
    \begin{figure}[!h]
    \begin{minipage}{.5\textwidth}
  \centering
  \includegraphics[width=0.8\textwidth]{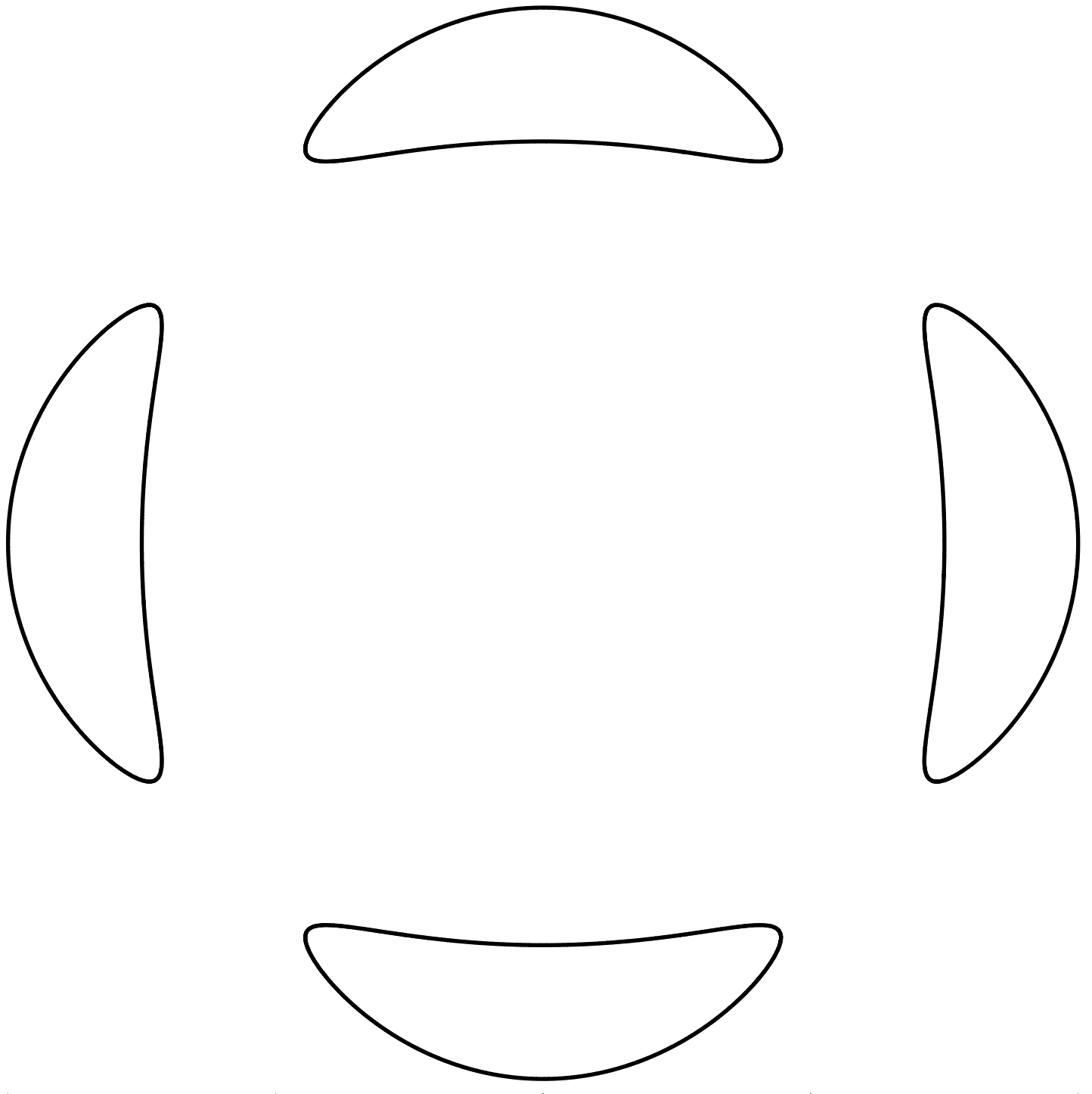}
\end{minipage}%
\begin{minipage}{.5\textwidth}
  \centering
  \includegraphics[width=0.8\textwidth]{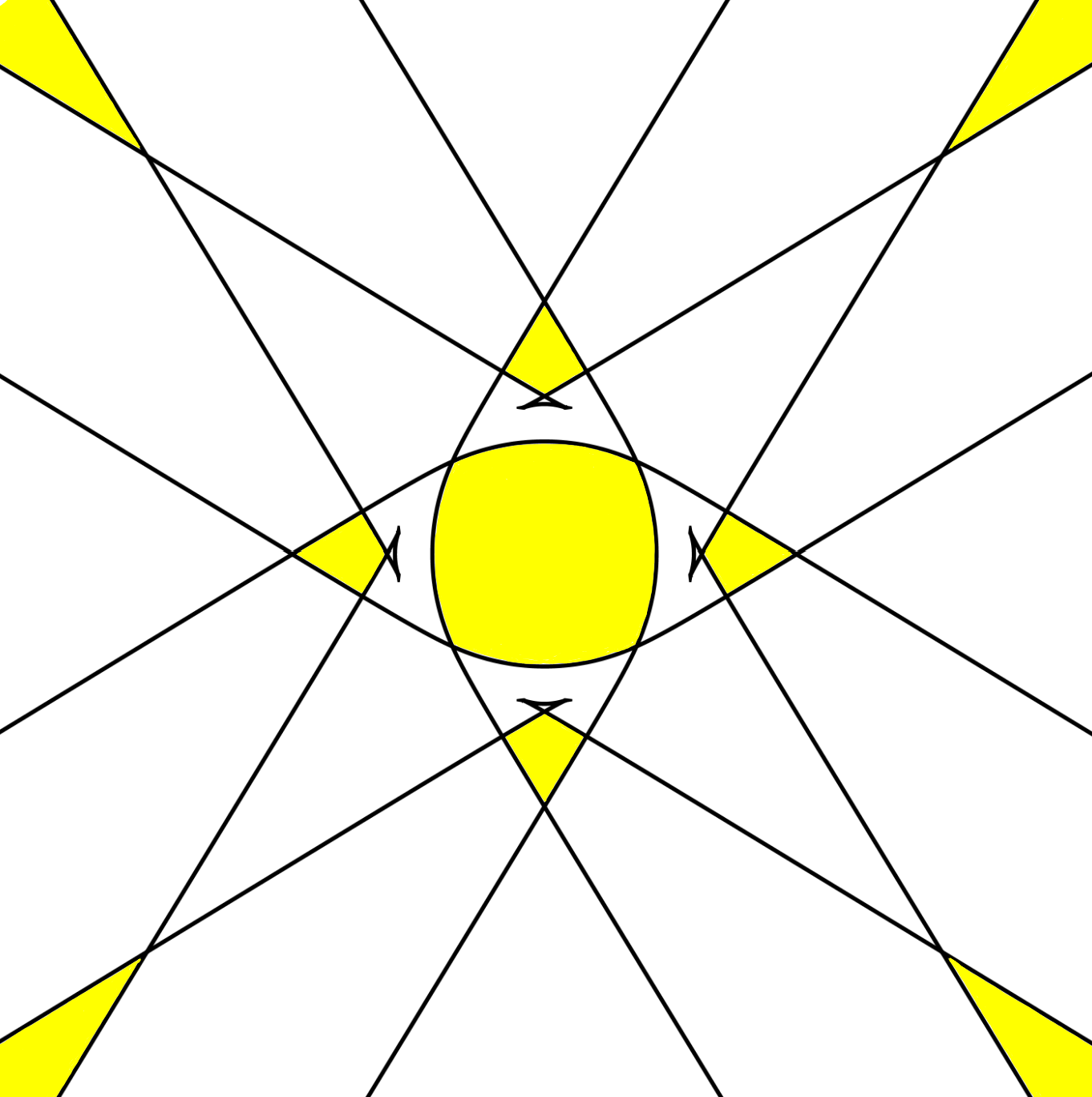}
\end{minipage}
        \caption{The Trott curve (left) and its avoidance locus (right)}\label{fig:trott}
    \end{figure}
\end{example}

A particularly important case of $\Acal_{n-1}(X)$ is when $X$ is the degree $2\ell$ Veronese embedding of $\PP^r$ in $\PP^{n-1},$ where $n = \binom{r+2k}{2k}.$ Then the avoidance locus $\Acal_{n-1}(X)$ is equal to the open cone of positive polynomials of degree $2\ell$, which has been studied extensively in the context of polynomial optimization~\cite[Chapter 3]{blekherman2012semidefinite}.
To see this, we first fix the notation $x_I = x_1^{i_1} \ldots x_r^{i_r}$ for a tuple $I = (i_1, \, \ldots, \,i_r).$ 
Then a hyperplane in $\PP^{n-1}$ intersecting $X$ may be pulled back along this embedding to the vanishing locus of $f(x) := \sum_{|I| = 2\ell}a_I x^I$. 
This hyperplane intersects $X$ in no real point if and only if the polynomial $f(x)$ is nonzero on all of $\RR^{r}.$

The structure of the paper is as follows. Section \ref{sec:properties} establishes basic properties of the avoidance locus. In Subsection \ref{subsec:discriminant}, we review the theory of higher Chow forms. Our main result in Section \ref{sec:properties} is Theorem \ref{thm: avoidance locus union of regions}, which states that the avoidance locus is a union of open regions in the complement of a higher Chow form. Subsection \ref{subsec:union_of_regions} builds up to the proof.

In Section \ref{sec:non adjacency}, we prove our first main result about avoidance loci: the Euclidean closures of the open regions making them up are not adjacent via any Zariski dense subset of the relevant higher Chow form. This can be seen concretely on the right in Figure \ref{fig:trott}, where the shaded regions are not adjacent.

Section \ref{sec:examples} focuses on computation and examples, including linear spaces, curves, and surfaces. In Section \ref{sec:n-1}, we focus on the case $k = n-1,$ in which case the avoidance locus lives in the dual projective space $(\PP^{n-1})^\vee$ and is equal to a union of regions in the complement of the dual variety $X^\vee$. We provide bounds for the number of regions it may contain in terms of the number of regions of $X_\RR.$ Code is available at \url{https://github.com/QWE123665/avoidance_locus}.

Finally, in Section \ref{sec:convexity} we explore convexity properties of the avoidance locus. We show in Theorem \ref{thm: hyperplane and X} that for the case $\cA_{n-1}(X),$ in which case the ambient Grassmannian is projective space, the avoidance locus is a union of \emph{convex} regions. For avoidance loci in more general Grassmannians, we study a property called \emph{slice convexity}. This is a notion of convexity first defined by Shamovich and Vinnikov to study hyperbolicity regions of algebraic varieties \cite{Shamovich_vinnikov}. Our second main result on avoidance loci is Theorem \ref{thm:slice convex}, which proves that that avoidance loci are unions of slice convex regions.

\section{Properties of the avoidance locus}\label{sec:properties}

In this section, we show that the avoidance locus of a variety is a union of regions in the complement of its corresponding discriminantal hypersurface in the Grassmannian. 
In particular, the avoidance locus is an open semi-algebraic subset of $\Gr(k,n)_\RR.$ We provide an introduction to Chow forms, Hurwitz forms and higher Chow forms in Subsection \ref{subsec:discriminant}. We prove  Subsection \ref{subsec:union_of_regions} that the avoidance locus is open and consists of regions in the complement of a higher Chow form. We begin by defining the following subvarieties of the Grassmannian which we will use throughout the paper.

\begin{definition}[Schubert varieties and lines in the Grassmannian]\label{def:schubert}
Let $\Gr(k,n)$ be the Grassmannian of projective $(k-1)$-dimensional linear spaces in $\mathbb{P}^{n-1}$. 
The following are linear spaces in the Plücker embedding of the Grassmannian, namely, cut out by linear forms in the Plücker variables. 
\begin{enumerate}
  \item[(a)] Fix a $(k-2)$-dimensional linear space $W \subset \mathbb{P}^{n-1}$.  
  The set
  \[
  \Gr(W,k) \; :=\; \{\, V \in \Gr(k,n) \mid W \subseteq V \,\}
  \]
  of $(k-1)$-spaces containing $W$ is a Schubert variety, linearly isomorphic to
  \(
   \mathbb{P}^{\,n-k}.
  \)

  \item[(b)] Fix a $k$-dimensional linear space $U \subset \PP^{n-1}$.  
  The set
  \[
  \Gr(k,U) \; :=\; \{\, V \in \Gr(k,n) \mid V \subseteq U \,\}
  \]
  of $(k-1)$-spaces contained in $U$ is a Schubert variety, linearly isomorphic to
  \(
  \mathbb{P}^k.
  \)

  \item[(c)] Fix $W \subset U$ with projective dimensions $k-2$ and $k.$ 
  The set
  \[
  L(W,U) \; :=\; \{\, [V] \in \Gr(k,n) \mid W \subseteq V \subseteq U \,\}
  \]
  is a \emph{line in $\Gr(k,n)$}, linearly isomorphic to $\bP^1.$
\end{enumerate}
\end{definition}

\subsection{Chow form, Hurwitz form and higher Chow forms}\label{subsec:discriminant}

Let $X \subset \PP^{n-1}$ be an irreducible projective variety of dimension $d$. Following \cite{Kohn}, we let $\text{Reg}(X)$ denote the set of smooth points of $X$ and $T_pX$ denote the embedded tangent space at the point $p$. 
\begin{definition}[Coisotropic variety]
For $i \in \{0, \, \ldots, \, d\},$ consider the incidence correspondence 
\[\Phi_i = \{(p, L) \, : \, p \in \text{Reg}(X \cap L), L\in \Gr(n-d+i-1,n) \text{ and } \dim (L \cap T_pX) \geq i)\}.\]
The \emph{$i$th coisotropic variety $\CH_i(X)$} of $X$ is Zariski closure of the projection of $\Phi_i$ to~$\Gr(n-d+i-1,n).$
\end{definition}
In words, each coisotropic variety consists of linear spaces which intersect $X$ non-transversely at a smooth point of $X.$ The higher Chow form $\CH_d(X)$ equals the dual variety $X^\vee,$ which parameterizes hyperplanes tangent to $X.$ The following corollary tells us that most coisotropic varieties are hypersurfaces.

\begin{proposition}[Corollary 6, Theorem 9 in \cite{Kohn}]\label{prop:hypersurface}
  The variety $\CH_i(X)$ has codimension one in $\Gr(n - d + i - 1, n)$ if and only if $0 \leq i \leq d - codim X^\vee + 1.$  Moreover, if $\CH_i(X)$ has codimension one, its degree equals the $i$-th polar degree of $X$.
\end{proposition}

When $i\geq 1$ and $\CH_i(X)$ is a hypersurface, we call its defining equation a \emph{higher Chow form}. The polar degrees of $X$, and hence the degrees of higher Chow forms, can be computed from the Chern classes of its tangent bundle via equation (3) in \cite{holme}. For certain values of $i$, higher Chow forms specialize to known discriminants. For example:

\begin{itemize}
    \item When $i=0,$ the variety $\CH_0(X)$ is always a hypersurface. 
It characterizes linear spaces of codimension $d+1$ in $\PP^{n-1}$ which meet $X.$ 
Its defining equation is called the \emph{Chow form} of $X$ and denoted $\Ccal_X.$ The degree of the Chow form equals the degree of $X$ itself. 
\item When $i = 1,$ the variety $\CH_1(X)$ is a hypersurface when $\deg X \geq 2$. Its defining equation, the \emph{Hurwitz form} $\Hcal_X$, characterizes linear spaces of complementary dimension that intersect $X$ at a non-reduced point. 
Let $g$ denote the sectional genus of $X$, i.e. the arithmetic genus of the curve $X\cap L$ where $L$ is a generic linear space of codimension $d-1$. %Then $\CH_1(X)$ is an irreducible hypersurface in $\Gr(n-d,n)$ defined by $\Hcal_X$. 
If the singular locus of $X$ has codimension at least two, then the degree of $\Hcal_X$ in Pl\"ucker coordinates is $2\deg X + 2g-2$ \cite[Theorem 1.1]{sturmfels2017hurwitz}.
\item As mentioned above, when $i=d$ we recover the definition of the dual variety $X^\vee$. A general hypersurface in projective space of degree $d$ has a dual variety of degree $d(d-1)^{n-1}.$
\end{itemize}

In the case of curves, the Hurwitz locus equals the dual variety, which parameterizes hyperplanes tangent to the curve. The case of curves in $\bP^2$ is particularly nice. A general plane curve $X$ of degree $d$ will have genus $\frac{(d-1)(d-2)}{2}.$ Thus $\deg \Hcal_X = 2d + (d-1)(d-2) -2 = d(d-1)$, which agrees with the formula for the degree of the dual curve.

\begin{example}[Rational normal curve]
Let $X$ be the rational normal curve, namely the image of $t\mapsto [1:t: \, \cdots \, : t^{n-1}]$ in $\PP^{n-1}$. Then a hyperplane section of $X$ corresponds to a univariate polynomial with degree $n-1$ and $\Hcal_X$ is the discriminant of univariate polynomials of degree $n-1$. 
The Hurwitz form has degree $2n-2$ since the genus of a rational normal curve is zero.
\end{example}

\begin{example}[Tact invariant] \label{eg:tact}
Suppose $X$ is the Veronese surface of $\PP^2$ in $\PP^5$, embedded via the standard embedding $[x:y:z] \mapsto [x^2: xy: xz: y^2: yz: z^2]$. 
The Hurwitz form $\Hcal_X \subseteq \Gr(4,6)$ encodes coefficients of two plane conics tangent at a point. The degree of the Veronese surface is $4$ and its sectional genus is 0 since a generic conic has genus 0, so $\deg \Hcal_X = 2\times 4-2 = 6.$
Its equation in dual Pl\"ucker coordinates is a degree 6 polynomial with 251 terms. If one parameterizes elements of $\Gr(4,6)$ as kernels of $2 \times 6$ matrices, then the tact invariant is the sum of $3210$ monomials in the twelve matrix entries \cite{breidingnotices}.
\end{example}

\begin{example}[Dual variety of Veronese surface]
Let $X$ be the Veronese surface as in Example \ref{eg:tact}.
The higher Chow form $\CH_2(X) \subseteq \Gr(4,6)$ encodes coefficients of hyperplanes tangent to $X$ and its corresponding hypersurface is the dual variety of $X$. 
The variety $\CH_2(X)$ has degree $3$ in the six Pl\"ucker coordinates of $\PP^5$ and its equation is 
\begin{equation}\label{eq:symdet}
    p_2^2p_3-p_1p_2p_4+p_0p_4^2+p_1^2p_5-4p_0p_3p_5.
\end{equation}
The Veronese embedding may also be viewed as the variety of rank one symmetric matrices inside the five-dimensional projective space of $3 \times 3$ symmetric matrices, with coordinates
\[\begin{bmatrix}
    2p_0 & p_1 & p_2 \\
    p_1 & 2p_3 & p_4 \\
    p_2 & p_4 & 2p_5
\end{bmatrix}.\]
The Hurwitz form in Equation \eqref{eq:symdet} is exactly the determinant of this matrix, up to scale.
\end{example}

\begin{example}[Segre variety in $\PP^3$]\label{ex: segre}
Suppose $X$ is the Segre variety of $\PP^1 \times \PP^1$ in $\PP^3$, embedded via the standard embedding $[x_0: x_1] \times [y_0: y_1] \mapsto [x_0y_0: x_0y_1: x_1y_0: x_1y_1]$. 
% The Hurwitz locus $\Hcal_X \subseteq \Gr(2,4)$ encodes coefficients of curves with support = the square.
We have $\deg \Hcal_X = 2 \times 2-2 + 0 = 2$ since a hyperbola has genus 0 and $\deg X = {2 \choose 1} = 2.$ Its equation in Pl\"ucker coordinates is
\begin{equation}\label{eq:hyperdet}
    p_{12}^2+p_{03}^2-2p_{02}p_{13}-2p_{01}p_{23}.
\end{equation}
Write a point in $\Gr(2,4)$ as the kernel of a $2 \times 4$ matrix. Then expression \eqref{eq:hyperdet} expanded in the eight matrix entries is the \emph{hyperdeterminant} of a $2 \times 2 \times 2$ tensor \cite[Equation 1.5]{gkzhyperdet}.
\end{example}

\begin{example}[Torus in $\PP^3$]\label{eg:torushurwitz}
We consider the surface $X$ in $\PP^3$ defined by 
\begin{align*}
F = (x^2+ y^2 + z^2 + 3 w^2)^2- 16(x^2+y^2)w^2.
\end{align*}
Its real part is topologically a torus.
The singular locus is a curve, so the degree formula for the Hurwitz form does not apply. However, a calculation shows that $\Hcal_X$ has degree $8$ and contains $111$ monomials when written in the standard basis of Pl\"ucker monomials on $\Gr(2,4).$
\end{example}

The following lemma relates the higher Chow forms of a projective variety to the higher Chow forms of a general linear section. We will need it in Section \ref{sec:non adjacency}.

\begin{lemma}[Slicing higher Chow loci]\label{lem:chowslice}
    Let $X \subset \bP^{n-1}$ be a $d$-dimensional projective variety. Let $\CH_i(X)$ be a higher Chow form of $X$ for some $0 \leq i \leq d,$ and set $k:= n-d-1+i.$ Suppose that $A$ is a $(k+\ell)$-dimensional linear subspace of $\bP^{n-1}$ intersecting $X$ transversally, where $\ell \geq 0.$ Then 
    \[CH_i(X) \cap \Gr(k,A) = \CH_{i}(X \cap A)\]
    as subvarieties of $\Gr(k,A) \cong \Gr(k,k+\ell+1).$
\end{lemma}

\begin{proof}
    Since the intersection of $A$ and $X$ is transverse, we have that $T_p(A \cap X) = T_pA \cap T_pX = A \cap T_pX.$ 
    We show that the complements in $\Gr(k,A)$ of each side of the equality are equal. That is, we claim that for $B \subset A$ of dimension $(k-1)$, the intersection $(X \cap A) \cap B$ is transverse in $A$ if and only if the intersection $X \cap B \subset A$ is transverse in $\bP^{n-1}.$ Let $A'$ and $B'$ denote the intersections of $A$ and $B$ with any affine chart of $\bP^{n-1}$ containing $p.$  
    In symbols, we are claiming that 
    \begin{equation}\label{eq:transverse}
        (A' \cap T_pX) + B' = A' \quad \iff \quad B' + T_pX = T_p\bP^{n-1}.
    \end{equation}
    Suppose the left side of \eqref{eq:transverse} holds. We take the direct sum with $T_pX$ and obtain that $T_pX + B' = T_pX + A',$ which is $T_p\bP^{n-1}$ by the transversality of $A \cap X.$ 

    Suppose the right side of \eqref{eq:transverse} holds. Then intersecting with $A',$ we obtain 
    $$(A' \cap B') + (A' \cap T_pX) = A' \cap T_p\bP^{n-1} = A'.$$ We conclude by recalling that $A'$ contains $B'.$
\end{proof}

\begin{example}[Slicing the torus]
    Consider the intersection of the torus in Example \ref{eg:torushurwitz} with the plane $A$ given by $x - 3y + 3z = 0.$ In the plane $A$, the intersection consists of an irreducible quartic curve. Its real component consists of two disjoint ovals, and its avoidance locus has two components. Intersecting the Hurwitz form of the torus in $\Gr(2,4)$ with the Schubert variety $\Gr(2,A)$ gives a degree-eight polynomial in the dual projective space $\Gr(2,A) \cong (\bP^2)^\vee.$ It is equal to the dual curve of the quartic curve. 
    \begin{figure}[!h]
    \centering
    \includegraphics[width = 0.4\textwidth]{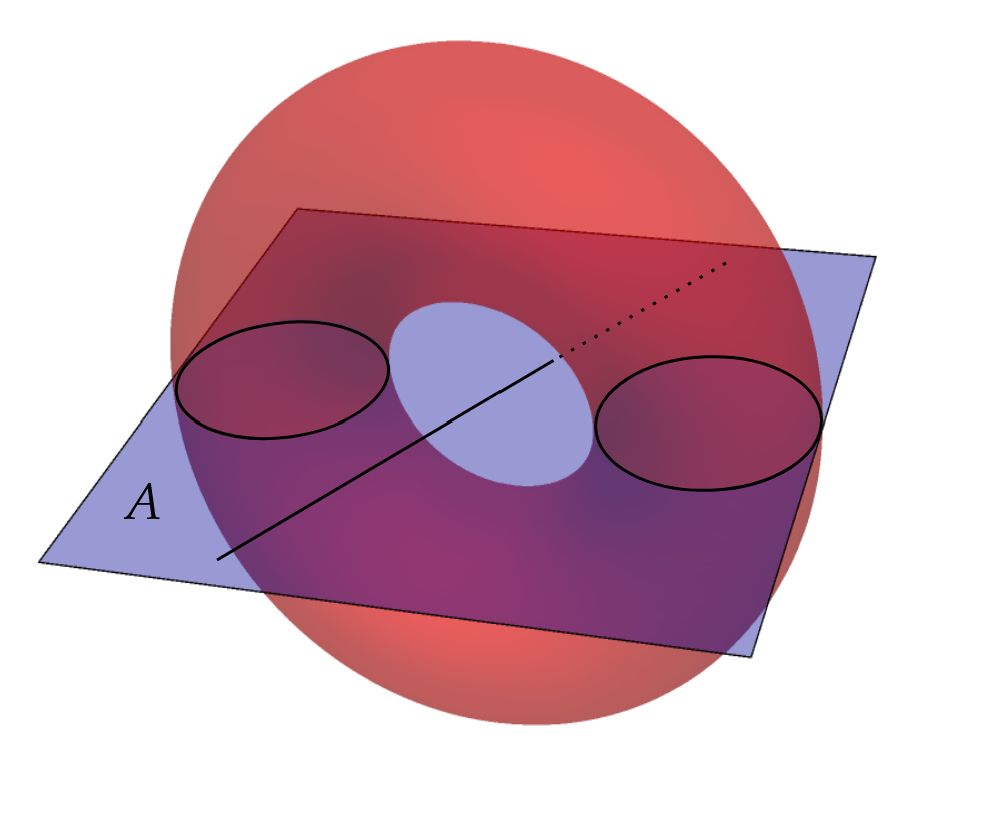}
    \caption{A line in the avoidance locus of the quartic curve obtained by slicing a torus}
    \end{figure}
    
    Lemma \ref{lem:chowslice} may fail if the linear space $A$ does not intersect $X$ transversely. Consider the intersection of the torus with the plane $A$ given by $x = 0.$ The result is two conics, with equations $3y^2-z^2-4yw+w^2=0$ and $3y^2-z^2+4yw+w^2=0.$ Intersecting the Hurwitz form in $\Gr(2,4)$ with the Schubert variety $\Gr(2,A)$ gives a degree eight polynomial in dual projective space $\Gr(2,A) \cong (\bP^2)^\vee.$ It factors as the product of the duals of the two conics, which each has degree two, and a term $(p_{12}^2 + p_{13}^2)^2.$
\end{example}

\subsection{The avoidance locus is open}\label{subsec:union_of_regions}
In this section, we consider the avoidance locus $\Acal_k(X)$ for $k\geq n-d-1$. 
Our main result is Theorem \ref{thm: avoidance locus union of regions}, which proves that the avoidance locus is a union of open regions in the complement of a higher Chow form. We start by establishing several basic properties of the avoidance locus.
\begin{lemma}\label{lem:easy_properties}
    Let $X \subset \PP^{n-1}$ be a smooth non-degenerate totally real variety of dimension $d$. The following properties of $\Acal_k(X)$ hold.
    \begin{enumerate}
        \item If the degree of $X$ is odd, then $\Acal_k(X) = \varnothing$ for $k \geq n-d$.
        \item If $\Acal_k(X)$ is nonempty, then $\Acal_{k-1}(X)$ is nonempty.
        \item (Union) If $Y$ is another non-degenerate smooth projective variety, then \[\Acal_k(X \cup Y) = \Acal_k(X) \cap \Acal_k(Y).\]
        \item (Intersection) If $Y$ is another non-degenerate smooth totally real variety, then \[\Acal_k(X) \cup \Acal_k(Y) \subseteq \Acal_k(X \cap Y).\]
        \item (Inclusion-reversing) If $Y \subset X$, then $\Acal_k(X) \subset \Acal_k(Y).$
    \end{enumerate}
\end{lemma}
\begin{proof}
For property (1), suppose that $L \cap X$ is purely complex for some linear space $L$, and that $L$ intersects $X$ transversely. Then its points all come in conjugate pairs. Intersecting $L \cap X$ with a general linear space of codimension equal to $\dim L \cap X$ produces a collection of complex conjugate points, which must be of even size. Since $L$ intersects $X$ transversely, this set of points has size equal to the degree $X.$ We leave properties (2)-(5) to the reader.
\end{proof}

\begin{proposition} 
For $k = n-d-1$, the avoidance locus $\Acal_k(X)$ is $\Gr(k,n)_\RR \setminus \CH_0(X)$.
\end{proposition}
\begin{proof}
The set $\Gr(k,n)_\RR \setminus \CH_0(X)$ is the collection of real linear spaces of dimension $k-1$ not intersecting $X$. So, we immediately obtain $\Acal_k(X) \supseteq \Gr(k,n)_\RR \setminus \CH_0(X)$.
The hypersurface $\CH_0(X)$ describes $(k-1)$-dimensional linear spaces intersecting $X$ non-transversely, and since $k-1 + \dim X$ is one less than the ambient space dimension, non-transversal intersection just means non-empty intersection. By definition of $\Acal_k(X)$, we have~$\Acal_k(X) \subseteq \Gr(k,n)_\RR \setminus \CH_0(X)$, which concludes the proof.
\end{proof}

For the rest of the section we build up to Theorem \ref{thm: avoidance locus union of regions}, which states that $\cA_k(X)$ is a union of open regions in the complement of a higher Chow form. We need the following definition of unavoidance locus to state our results.
Our strategy is to show that the unavoidance locus and the avoidance locus are both open in the Euclidean topology, which implies that the boundary of the avoidance locus is contained in the higher Chow form.

\begin{definition}
Let \( X \subset \mathbb{P}^{n-1} \) be a smooth non-degenerate totally real variety of dimension $d$. The \( k \)-th \emph{unavoidance locus} is
\[
\mathcal{U}_k(X) := \left\{ V \in \operatorname{Gr}(k,n)_\RR : V \cap X \text{ is transverse and } (V \cap X_\mathbb{R}) \neq \emptyset \right\}.
\]
\end{definition}

In particular, the real Grassmannian $\Gr(k,n)_\RR$ breaks up as a disjoint union
\begin{equation}\label{eq:sqcup}
    \Gr(k,n)_\RR \ = \ \Ucal_k(X) \, \sqcup \, CH_{k+d-n+1}(X)_\RR \, \sqcup \, \Acal_k(X).
\end{equation}

We review Bertini's Theorem \cite[Theorem 17.16]{harris2013algebraic}, which will be used throughout the rest of the paper.
It states that for a smooth variety, a general member of a linear system of divisors on $X$ is smooth away from the base locus of the system. 
It extends to the following, see \cite[Corollary 2.7]{ranestad2024real}.

\begin{corollary}[of Bertini's Theorem]\label{cor:extension of bertini}
Let $X$ be a smooth variety. 
If the base locus of a linear system of divisors on $X$ is empty or a finite reduced set of points, then the general member of the linear system is smooth. 
\end{corollary}

We will also need a lemma which generalizes the result that the roots of a polynomial are locally continuous functions in its coefficients.

\begin{lemma}\label{lem:points continuous}
Let $X \subset \bP^{n-1}$ be a projective variety and let $V$ be a general linear space of complementary dimension. Then \( V \cap X \) is a collection of $\deg X$ many distinct points whose coordinates are locally continuous functions of the Pl\"ucker coordinates of $V.$
\end{lemma}
\begin{proof}
Roots of a univariate polynomial change continuously as its coefficients change. 
By considering the Gr\"obner basis of the ideal $I(V \cap X)$ generated by polynomials defining $X$ and the linear relations defining $V$, the $i$-th coordinate of the solutions to $V\cap X$ are roots of some univariate polynomial with coefficients that change continuously as we move $V$ within a small neighborhood. 
Hence, the solutions to $V\cap X$ change continuously as we move $V$.
\end{proof}

\begin{lemma}\label{lem: unavoidance locus open}
Let \( X \) be a non-degenerate totally real smooth variety in \( \mathbb{P}^{n-1} \) of dimension $d$. 
For \( k \geq n-d \), the unavoidance locus \( \mathcal{U}_k(X) \) is non-empty and open in the Euclidean topology.
\end{lemma}
\begin{proof}
We start by proving that \( \mathcal{U}_k(X) \) is nonempty. By assumption, there exists a smooth real point~\( p \in X_\mathbb{R} \). Let \( V \) be a general linear space of projective dimension \( k-1 \) in \( \mathbb{P}^{n-1} \) containing $p.$ The base locus of the linear system of $(k-1)$-spaces containing $p$ is the point $p$ itself. Thus by Corollary \ref{cor:extension of bertini}, \( V \cap X \) is transverse. Since $p$ is in $ V\cap X$ we have \( (V \cap X)_\mathbb{R} \neq \emptyset \). Thus $[V]$ is in the unavoidance locus. % Hence, \( [V] \in \mathcal{U}_k(X) \).

We now prove that \( \mathcal{U}_k(X) \) is open. Suppose first that \( k = n-d \). 
If \( [V] \in \mathcal{U}_k(X) \), then \( V \cap X \) is a collection of $\deg X$ many distinct points with at least one of them real. 
By Lemma \ref{lem:points continuous}, the coordinates of the intersection points are locally continuous functions in the Pl\"ucker coordinates of $V.$ Thus there exists a small open neighborhood of $[V]$ such that if we slightly perturb \( [V] \) to \( [V'] \), \( V' \cap X \) is still a collection of $\deg X$ many distinct points and a real point in \( V \cap X \) will move to a real point in \( V' \cap X \). 
Thus the set \( \mathcal{U}_{k}(X) \) is open.

Now we suppose that \( k > n-d \) and \( [W] \in \mathcal{U}_k(X) \). We will reduce to the earlier case $k=n-d$. 
Let $p$ be a generic point in $(W\cap X)_\RR$.
Consider the linear system of \( [V] \in \operatorname{Gr}(n-d,n)_\RR \) with $p\in V \subset W.$ 
If $V$ is a general element of this linear system, then  $W \cap X$ is smooth by definition of the unavoidance locus. Thus \( V \cap X \) is smooth by the extension of Bertini's Theorem. 
We also have that \(  (V \cap X)_\mathbb{R} \) is nonempty, since it contains $p$.   

Let us identify \( W \) with the row span of a \( k \times n \) matrix \( M_W \) and \( V \) with the row span of a \( (n-d) \times n \) matrix \( M_V \). 
There exists a matrix \( M \in \mathbb{R}^{(n-d) \times k} \) such that \( M_V = M M_W \). The matrix \( M \) has full row rank and it defines a continuous map
\begin{align*}
    \varphi: \Gr(k,n)_\RR & \to \Gr(n-d,n)_\RR \\
    [A] & \mapsto [M\cdot A].
\end{align*}
In particular, \( \varphi([W]) = [V] \). 
We have shown that there exists an open neighborhood \( \mathcal{V} \) of \( [V] \) such that for all \( [V'] \in \mathcal{V} \), the intersection \( V' \cap X \) is transverse and \( (V' \cap X)_\mathbb{R} \neq \emptyset \). 
Thus,
\(\varphi^{-1}(\mathcal{V})\)
is a non-empty open set in \( \operatorname{Gr}(k,n)_\RR \) such that any linear space in \( \varphi^{-1}(\mathcal{V}) \) intersects \( X_\mathbb{R} \).
We conclude that
\(
\varphi^{-1}(\mathcal{V}) \cap \CH_{k+d-n+1}(X)^c 
\)
is an open neighborhood of \( [W] \) contained in \( \mathcal{U}_k(X) \), and hence \( \mathcal{U}_k(X) \) is open.
\end{proof}

\begin{lemma}\label{lem: avoidance locus open}
Suppose \( X \subset \bP^{n-1}\) is a smooth non-degenerate totally real variety of dimension $d$. 
For \( k \geq n-d \), the avoidance locus \( \mathcal{A}_k(X) \) is open.
\end{lemma}
\begin{proof}
If \( \Acal_{k}(X)\) is empty then the statement is vacuously true. Now, suppose that \( \Acal_{k}(X)\) is nonempty. We work in $\mathbb{S}^{n-1}$, the unit sphere in the Euclidean space $\RR^{n}$. We define $X_1, \dots, X_r$ to be the intersections with $\mathbb{S}^{n-1}$ of the affine cones in $\RR^{n}$ above the components of $X_\RR.$
Each $X_i$ is a smooth manifold of dimension equal to $\dim X$ \cite[Proposition 7.6.2]{bochnak2013real}. Since it is a closed subset of the unit sphere, it is compact.

We will show that for $[V]\in \Acal_k(X)$, there exist an open neighborhood $\Vcal$ of $[V]$ lying in $\Acal_k(X).$
We define on $\Gr(k,n)_\RR$ the function
\[
\phi_i(W) := \inf_{x \in X_i} \operatorname{dist}(x, W)
\]
where $\text{dist}$ is the geodesic distance in $\mathbb{S}^{n-1}$ from the point $x$ to the affine cone of the real linear space $W$ intersecting with $\mathbb{S}^{n-1}$. 
Since the distance function is continuous and each \( X_i \) is compact, the infimum is attained, so the function $\phi_i$ is well-defined. 
We have $\phi_i(V)>0$ for $i=1,\ldots,r$.
By continuity, there exists an open neighborhood $\mathcal{V} $ of $V$ in $\Gr(k,n)_\RR$  such that the function value of each $\phi_i(\cdot)$ will remain positive. 
It follows that every $[V]$ in $\Vcal$ avoids the variety $X_\RR$. Then $\Vcal\cap (\CH_{k+d-n+1})^c$ is open, it contains $[V]$ and it is contained in $\Acal_k(X)$, so it is the desired open neighborhood. 
\end{proof}

Lemmas \ref{lem: unavoidance locus open} and \ref{lem: avoidance locus open} together tell us that if $\CH_{i}(X)$ is not a hypersurface, then the avoidance locus $\cA_{n-d+i-1}$ is empty.  We can see this behavior for instance for Segre embeddings of products of projective spaces. Note that by Proposition \ref{prop:hypersurface}, if $\CH_i(X)$ is not a hypersurface for some $0 \leq i \leq d$, then the dual variety $X^\vee$ is also not a hypersurface.

\begin{example}
    The degree of the Segre embedding of $X:= \bP^{k_1} \times \ldots \times \bP^{k_r}$ is the multinomial coefficient $\frac{(k_1 + \ldots + k_r)!}{k_1! \ldots k_r!}.$ In particular, if this degree turns out to be odd, then all avoidance loci are empty. Let $k:= \sum_{i=1}^r k_i.$ The dual variety is a hypersurface if and only if $2k_i \leq k$ for each~$i$  \cite{GKZ}. Otherwise, the codimension of $X^\vee$ is $1 + 2k_i - k$  \cite{kaji}.

    For example, consider the Segre embedding $\PP^1 \times \PP^3 \hookrightarrow \PP^7.$ The image has even degree, namely four, but the dual variety has codimension three by the degree formula above. Thus $\mathcal{A}_6(X)$ and $\cA_7(X)$ are empty.
\end{example}

\begin{theorem}\label{thm: avoidance locus union of regions}
Suppose \( X \subset \bP^{n-1} \) is a smooth non-degenerate totally real variety of dimension $d$. 
For \( k \geq n-d \), the avoidance locus \( \mathcal{A}_k(X) \) is an open set in \( \Gr(k,n)_\RR \), consisting of the union of regions in
\[
\Gr(k,n)_\RR \setminus \CH_{k+d-n+1}(X).
\]
\end{theorem}
\begin{proof}
We have shown in Lemma \ref{lem: unavoidance locus open} and \ref{lem: avoidance locus open} that both the unavoidance locus $\Ucal_k(X)$ and the avoidance locus $\Acal_k(X)$ are open. 
The real Grassmannian decomposes as a disjoint union $\Gr(k,n)_\RR \ = \ \Ucal_k(X) \, \sqcup \, \CH_{k+d-n+1}(X)_\RR \, \sqcup \, \Acal_k(X)$, so the boundary of $\Acal_k(X)$ lies in $\Ucal_k(X)\sqcup \CH_{k+d-n+1}(X)$. 
Assume for contradiction that $[V]\in \partial \Acal_k(X)$ and $[V]\in \Ucal_k(X)$. Since $\Ucal_k(X)$ is open, there exists an open neighborhood $\Vcal$ of $[V]$ such that $\Vcal\subseteq \Ucal_k(X)$ and in particular $\Vcal \cap \Acal_k(X)$ is empty. This is a contradiction.
\end{proof}

\section{Non-adjacency of Regions}\label{sec:non adjacency}
In this section our main result is Theorem \ref{thm:nonadjacent}, which states that the regions in the avoidance locus are not adjacent via any Zariski dense subset of the relevant higher Chow form. We prove this theorem in Subsection \ref{subsec:not adjacent}. In Subsection \ref{subsec:single region}, we show that for a fixed variety  $X$, avoidance loci of higher-dimensional linear spaces give us information about avoidance loci of lower-dimensional linear spaces. 
We then prove in 
Proposition \ref{prop:singleregion} that the lower-dimensional linear spaces obtained from a region of $\cA_k(X)$ this way are contained in a single region of the avoidance locus $\cA_{k-1}(X).$

\subsection{Regions in $\mathcal{A}_k(X)$ are not adjacent}\label{subsec:not adjacent}

We call two regions in the complement of a higher Chow form $\CH_i(X)$ \emph{adjacent} if the intersection of their Euclidean closures is Zariski dense in $\CH_i(X)$.
In this subsection we prove that the avoidance locus $\Acal_k(X)$ consists of \emph{non-adjacent} regions in the complement of a higher Chow form, culminating in Theorem \ref{thm:nonadjacent}. 
The $k = n-d$ case was proved in \cite{ranestad2024real}. 
The key idea is to use the relationship between $\Acal_k(X)$ and $\Acal_{k-1}(X)$ to reduce the problem to the case of hyperplanes.

It is a result of Holweck that any bitangent hyperplane to a variety $X$ is a component of the singular locus of the dual variety $X^\vee$ \cite{holweck}. In Proposition \ref{prop:bitangent} we show that similarly, linear spaces tangent to $X$ at multiple points are singular points of the higher Chow locus. This fact will be used in the proof of Theorem \ref{thm:nonadjacent}.

\begin{proposition}\label{prop:bitangent}
Suppose \( X \subset \bP^{n-1}\) is a smooth non-degenerate totally real variety of dimension $d$. 
If $[V]$ is a smooth point on $\CH_{k+d-n+1}(X)$ where $k \geq n-d$, then $V$ intersects $X$ non-transversely at exactly one point.
\end{proposition}

\begin{proof}
We identify $\bP(\CC^{n-k} \otimes \CC^k)$ with the (projective) Stiefel manifold of rank $n-k$ matrices of size $(n-k)\times n$. Inside this Stiefel manifold we have $\PP^{n-k-1} \times \PP^{n-1} \hookrightarrow \PP(\CC^{n-k} \otimes \CC^n).$

The higher Chow form $\CH_{k+d-n+1}(X)$ can be identified with $ (\mathbb{P}^{n-k-1} \times X)^\vee$  via the map
\[
\pi : St(n-k, n) \longrightarrow Gr(k, n), \quad A \longmapsto \PP(\ker A).
\]
Indeed, $A \in (\mathbb{P}^{n-k-1} \times X)^\vee$ is tangent to $\mathbb{P}^{n-k-1} \times X$ at some point $(y,x)$ if and only if $\PP( \ker A)$ intersects $X$ non-transversally at $x$, by \cite[Proposition 5]{Kohn}.

It suffices to show the result for hyperplanes intersecting the smooth variety $\mathbb{P}^{n-k-1} \times X$ inside $\PP(\CC^{n-k} \otimes \CC^n)$. By the biduality theorem (\cite[Theorem 1.1]{GKZ}), if a hyperplane $H$ is tangent to $\mathbb{P}^{n-k-1} \times X$ at $p$, then the hyperplane $H_p$ in $(\mathbb{P}^{n-k-1} \times X)^\vee$ corresponding to $p$ is tangent to the point $[H]$. Since $(\mathbb{P}^{n-k-1} \times X)^\vee$ is a hypersurface, through a smooth point $[H]$ there is a single tangent hyperplane $H_p$. Thus $H$ is tangent to $\mathbb{P}^{n-k-1} \times X$ only at $p.$
\end{proof}

We now prove the non-adjacency result for the avoidance locus of hyperplanes. This can be seen concretely for the avoidance locus of the Trott curve in Figure \ref{fig:trott}, in which no two yellow regions share a border. We note that the assumption that $X$ is positive-dimensional is necessary; otherwise $X^\vee$ will be a union of linear spaces, and the avoidance locus of hyperplanes will be equal to $(\bP^{n-1})^\vee \setminus X^\vee.$

\begin{proposition}\label{prop:nonadjacent hyperplane}
    Suppose \( X \subset \bP^{n-1} \) is a smooth non-degenerate totally real variety of dimension $d$.  Two distinct regions in $\cA_{n-1}(X)$ are not adjacent, i.e. the intersection of their Euclidean closures is not Zariski-dense in $X^\vee.$
\end{proposition}

\begin{proof}
Suppose that $\cA_{n-1}(X)$ is non-empty. Let $X_1, \, \ldots, \, X_r$ be the connected components of $X_{\RR}$. 
As in the proof of Proposition \ref{lem: avoidance locus open}, we define the function
\[
\phi_i(H) = \inf_{x \in X_i} \mathrm{dist}(x,H)
\]
for hyperplanes $H$ in $\bP^{n-1}$. Suppose we have two regions $\mathcal{V}_1$ and $\mathcal{V}_2$ of the avoidance locus which share a generic point $[V]$ of $X^\vee.$ Then any open neighborhood of $[V]$ will intersect $\Vcal_1, \Vcal_2$ in open sets. Thus there exist points $H_1 \in \mathcal{V}_1, H_2 \in \mathcal{V}_2$ and a curve $\gamma$ passing through $[V]$ which connects $H_1$ and $H_2$ and does not pass through any other point of $X^\vee.$ Then $\phi(\gamma(t))$ is positive except for some value $t_*$ for which $\gamma(t_*) = [V].$  

If $H_1$ and $H_2$ are in different components of the avoidance locus, then there must be some component $X_i$ for which the signs of the linear forms $\ell_1$ and $\ell_2$ defining $H_1$ and $H_2$ differ. 
Suppose without loss of generality that this is true for $X_1.$ Then any $x$ in the interior of $\text{conv}(X_1)$ satisfies $\ell_1(x) >0$ and $\ell_2(x) < 0.$ If we fix $x,$ then by continuity there must be some $s$ such that the linear form defining $\gamma(s)$ is zero at $x.$ Then the hyperplane $\gamma(s)$ meets $X_1,$ so $s = t_*$ and $\gamma(t_*)$ contains $x$. Since the same argument works for any $x$ in the interior of $X_1,$
the hyperplane $\gamma(t_*)$ also contains $X_1$. 
But $\gamma(t*)$ is a generic point on $X^\vee$, so by Proposition \ref{prop:bitangent}, it can only be tangent to $X$ at a point. This implies that $X_1$ is a point. But $X_1$ has dimension $d>0$, which is a contradiction.
\end{proof}

We now extend the non-adjacency result from Proposition~\ref{prop:nonadjacent hyperplane} to $\cA_k(X)$ for $n-d-2 <k< n-1,$ namely to avoidance loci other than the complement of the Chow form.

\begin{theorem}\label{thm:nonadjacent}
Suppose \( X \subset \bP^{n-1}\) is a smooth non-degenerate totally real variety of dimension $d$. 
Fix $n-d-2 < k \leq n-1.$ Two distinct regions in $\mathcal{A}_k(X)$ are not adjacent via any Zariski dense subset of $\CH_{k+d-n+1}(X)$.
\end{theorem}

\begin{proof}
Suppose that $k<n-1$ and that $\mathcal{V}_1, \mathcal{V}_2$ are two regions in $\mathcal{A}_k(X)$ whose boundaries share a real smooth point $[V]$ of $\CH_{k+d-n+1}(X)$. We choose a general $k$-space $A$ containing $V$ such that the Schubert variety $\Gr(k,A)$ satisfies
\begin{enumerate}[label=\roman*.]
    \item $\Gr(k,A)_\RR$ is not contained in $\mathcal{V}_1$ or $\mathcal{V}_2$ and $\mathcal{V}_1\cap \Acal, \mathcal{V}_2 \cap \Acal\neq \emptyset$
    \item $\Gr(k,A)$ is not contained in $\CH_{k+d-n+1}(X)$, 
    \item $[V]$ is a smooth point of $Y:= \Gr(k,A) \cap \CH_{k+d-n+1}(X)$, and
    \item $A \cap X$ is smooth. 
\end{enumerate}

Condition (i) is an open condition in the Euclidean topology. We claim that the other conditions are Zariski open, beginning with (ii). By Lemma \ref{lem:chowslice}, we have that $\Gr(k,A) \cap \CH_{k+d-n+1}(X) = \CH_{k+d-n+1}(A \cap X).$ 
This is all of $\Gr(k,A)$ if and only if $A \cap X = A,$ i.e. $A\subset X$. But containment within the variety $X$ is a Zariski closed condition on $k$-spaces. 

Condition (iii) is Zariski open by Bertini's theorem. There is a map from $\RR^n/V$ to linear spaces in $\Gr(k,n)_\RR$ which sends $v$ to $\Gr(k, \bP(V \oplus \langle v \rangle)).$ This map is linear in the entries in $v,$ so it indeed sweeps out a linear system on $\CH_{k+d-n+1}(X)$. The base locus of this system is exactly $[V],$ so a general linear section in it is smooth by Corollary \ref{cor:extension of bertini}, and in particular smooth at $[V]$. Each possible $A$ is of the form $V \oplus \langle v \rangle),$ so this gives us a Zariski open set of $A$ satisfying condition (iii).

Finally we argue (iv), namely that if $A$ is Zariski general, then $A \cap X$ is smooth. 
Since $A$ is a general linear space containing $V$, by Bertini's Theorem, we only need to show that $A\cap X$ is smooth at $\Sing(V\cap X)$.
Since $V$ is a smooth point of the higher Chow locus, the intersection $V\cap X$ is singular at one point $p$ by Proposition \ref{prop:bitangent}. We have $\dim(V \cap T_pX) = 1+ \dim X- \codim V$, which implies that $V\oplus T_p X$ is not all of $\PP^{n-1}$.  
Thus we can choose $A$ such that $A = V \oplus \langle q \rangle $ where $q\notin V\oplus T_p X,$ then $\dim(A \cap T_pX) = 1+ \dim X- \codim V = \dim X- \codim A$ and $A \cap X$ is smooth at $p$.

We have three Zariski open conditions and one Euclidean open condition, and thus we may indeed make a choice of such $A.$ Now, we reduce to the dual variety case. By Lemma \ref{lem:chowslice}, we have that $(A \cap X)^\vee \subset A^\vee$ is contained in $Y.$ Thus $\mathcal{V}_1 \cap A$ and $\mathcal{V}_2 \cap A$ are adjacent via a smooth point $[V]$ of $A \cap X.$ Then $\mathcal{V}_1 \cap A$ and $\mathcal{V}_2 \cap A$ are each contained in the avoidance locus $\cA_k(A \cap X),$ which contradicts Proposition \ref{prop:nonadjacent hyperplane}. Note that $A \cap X$ has dimension $k - (n-d-2),$ which is greater than zero by assumption.
\end{proof}

\subsection{Relations between avoidance loci} \label{subsec:single region}

The avoidance locus $\cA_{k}(X)$ gives us information about avoidance loci of linear spaces of lower dimension. 
The idea is that if $V$ is in the avoidance locus $\cA_{k}(X),$ then any generic codimension one subspace of $V$ must be in the avoidance locus $\Acal_{k-1}(X).$ Essentially, we can get a subset of lower-dimensional subspaces avoiding $X$ for free. In Proposition \ref{prop:singleregion}, we show that each region of $\Acal_k(X)$ yields a single region of $\Acal_{k-1}(X)$ via this method. 

\begin{lemma}\label{lem: two grassmannians}
We equip each real Grassmannian with the Euclidean topology induced from its ambient projective space. Suppose $\mathcal{U}$ is an open set in $\operatorname{Gr}(k{+}1,n)_\RR$.
Define $\mathcal{U}_p$ as
\[
\mathcal{U}_p := \{\, [W] \in \Gr(k,n)_\RR \, \mid \, W \subseteq V \text{ for some } [V] \in \mathcal{U} \,\}.
\]
Then $\mathcal{U}_p$ is an open set in $\operatorname{Gr}(k,n)_\RR$. If $\mathcal{U}$ is connected, $\mathcal{U}_p$ is connected.
\end{lemma}

\begin{proof}
    Consider the diagram in Figure \ref{fig:bundle}, where $\rm{Fl}(k,k+1,n)$ is the partial flag variety parameterizing two-step flags $V \subset W$ of dimensions $k$ and $k+1$. 
    \begin{figure}[!h]
\begin{center}
\begin{tikzcd}
& \text{Fl}(k,k+1,n) \arrow{dl}[swap]{\pi_1} \arrow{dr}{\pi_2} & \\
\Gr(k,n)_\RR & & \Gr(k+1,n)_\RR 
\end{tikzcd}
\end{center}
\caption{Flags of two linear subspaces}\label{fig:bundle}
\end{figure}

The maps $\pi_1, \pi_2$ give $\text{Fl}(k,k+1,n)$ the structure of a smooth fiber bundle with fibers $\bP^k$ and $\bP^{n-k-1},$ respectively. They are open quotient maps  \cite[Exercise 10-19]{Lee}. Thus $\mathcal{U}_p = \pi_1(\pi_2^{-1}(\mathcal{U}))$ is open. Furthermore, the fibers are connected, so $\pi_2^{-1}(\mathcal{U})$ is connected.
\end{proof}

\begin{proposition}\label{prop:singleregion}
Suppose \( X \subset \bP^{n-1}\) is a smooth non-degenerate totally real variety of dimension $d$. 
If $\Ucal$ is a region of $\Acal_k(X)$, then the set
\[
\Ucal' := \{\, [W] \in \Gr(k-1,n)_\RR \mid W \subseteq U \text{ for some } [U] \in \Ucal,\, W\cap X \text{ transverse}\,\} \
\]
is contained in a single region of $\Acal_{k-1}(X)$.
\end{proposition}

\begin{proof}
The set $\mathcal{U}_p$ of $(k-2)$-dimensional linear spaces contained in some linear space in $\mathcal{U}$ is open and connected, by Lemma \ref{lem: two grassmannians}. 
Thus to show that $\mathcal{U}_p$ is contained in one region, it suffices to show it does not intersect a Zariski dense subset of the higher Chow form in $\Gr(k-1,n).$
This follows immediately from Theorem \ref{thm:nonadjacent}, otherwise there are two regions in $\Acal_{k-1}(X)$ that are adjacent via the higher Chow locus.
\end{proof}

\begin{figure}[!h]
    \centering
    \tikzset{every picture/.style={line width=0.75pt}} %set default line width to 0.75pt        

\begin{tikzpicture}[x=0.75pt,y=0.75pt,yscale=-1,xscale=1]

%Shape: Ellipse [id:dp579627456109501] 
\draw   (160,220) .. controls (160,197.91) and (177.91,180) .. (200,180) .. controls (222.09,180) and (240,197.91) .. (240,220) .. controls (240,242.09) and (222.09,260) .. (200,260) .. controls (177.91,260) and (160,242.09) .. (160,220) -- cycle ;
%Curve Lines [id:da4009315791808017] 
\draw    (160,220) .. controls (169.25,238.25) and (231.25,236.75) .. (240,220) ;
%Curve Lines [id:da24873239877177256] 
\draw  [dash pattern={on 0.84pt off 2.51pt}]  (160,220) .. controls (163.25,205.25) and (236.75,204.25) .. (240,220) ;
%Shape: Ellipse [id:dp8182754019296866] 
\draw   (260,80) .. controls (260,57.91) and (277.91,40) .. (300,40) .. controls (322.09,40) and (340,57.91) .. (340,80) .. controls (340,102.09) and (322.09,120) .. (300,120) .. controls (277.91,120) and (260,102.09) .. (260,80) -- cycle ;
%Curve Lines [id:da5858326485928973] 
\draw    (260,80) .. controls (269.25,98.25) and (331.25,96.75) .. (340,80) ;
%Curve Lines [id:da3526368830437815] 
\draw  [dash pattern={on 0.84pt off 2.51pt}]  (260,80) .. controls (263.25,65.25) and (336.75,64.25) .. (340,80) ;
%Shape: Ellipse [id:dp8029025497639176] 
\draw   (390.67,222) .. controls (390.67,199.91) and (408.58,182) .. (430.67,182) .. controls (452.76,182) and (470.67,199.91) .. (470.67,222) .. controls (470.67,244.09) and (452.76,262) .. (430.67,262) .. controls (408.58,262) and (390.67,244.09) .. (390.67,222) -- cycle ;
%Curve Lines [id:da15010632394679946] 
\draw    (390.67,222) .. controls (399.92,240.25) and (461.92,238.75) .. (470.67,222) ;
%Curve Lines [id:da5841069245128283] 
\draw  [dash pattern={on 0.84pt off 2.51pt}]  (390.67,222) .. controls (393.92,207.25) and (467.42,206.25) .. (470.67,222) ;
%Shape: Parallelogram [id:dp061436158742061764] 
\draw  [fill={rgb, 255:red, 221; green, 122; blue, 230 }  ,fill opacity=0.5 ] (265.9,166.78) -- (420,60) -- (402.4,175.67) -- (248.3,282.45) -- cycle ;
%Shape: Parallelogram [id:dp053246612600894694] 
\draw  [fill={rgb, 255:red, 126; green, 211; blue, 33 }  ,fill opacity=0.5 ] (221.85,164.08) -- (360,280) -- (343.77,179.1) -- (205.61,63.18) -- cycle ;
%Straight Lines [id:da5173160706070113] 
%\draw [color={rgb, 255:red, 208; green, 2; blue, 27 }  ,draw opacity=1 ]   (300.48,142.88) -- (311,239.6) ;
\draw [color={rgb, 255:red, 74; green, 144; blue, 226 }  ,draw opacity=1 ]   (300.48,142.88) -- (311,239.6) ;

\end{tikzpicture}
    \caption{An avoidant line (in blue) from avoidant hyperplanes in different regions}
    \label{fig:spheres}
\end{figure}
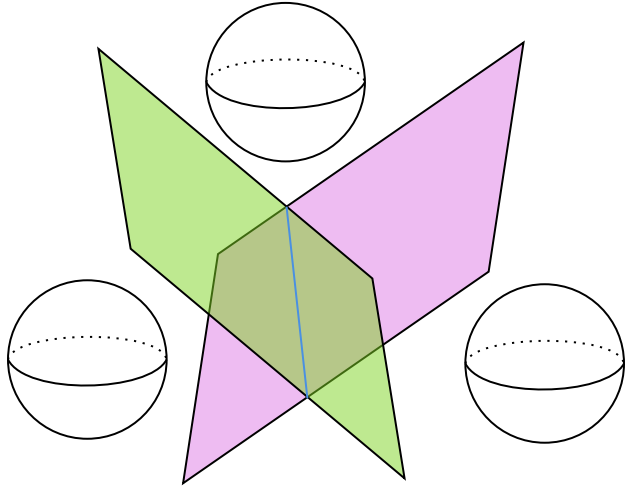

The reverse of Proposition \ref{prop:singleregion} is not true; one lower-dimensional linear space may come from two higher-dimensional linear spaces in different avoidance loci. An example of this is shown in Figure \ref{fig:spheres}, where $X$ is the union of three spheres. Thus, knowing the number of components in $\cA_{k}(X)$ does not necessarily give us bounds on the number of components in $\cA_{k-1}(X),$ other than telling us that the latter is nonempty.

\section{Computations and Algorithms}\label{sec:examples}

In this section, we illustrate how to compute avoidance loci in concrete examples, from linear spaces to curves and hypersurfaces.

\subsection{Linear spaces}

\begin{lemma} \label{lem:linear}
    Let $V $ be a $d$-dimensional linear space in $\PP^{n-1}.$ Then there are three cases for $\Acal_k(V),$ depending on $k.$
    \begin{enumerate}[i)]
        \item $k>n-d-1$: The avoidance locus is empty.
        \item $k < n-d-1$: The avoidance locus is the complement of a codimension $n-d-k$ Schubert variety in $\Gr(k,n)_\RR,$ and thus dense and connected.
        \item $k = n-d-1$: The avoidance locus is the complement of a Schubert hyperplane in $\Gr(n-d-1)_\RR.$
    \end{enumerate}
\end{lemma}

\begin{proof}
For (i), we note that any real linear space of dimension at least $n-d$ meets $V.$ For (ii) and (iii), the Schubert variety in question consists of $(k-1)$-dimensional linear subspaces of $\bP^{n-1}$ which meet $V$. It has codimension $n-d-k$ in $\Gr(k,n)_\RR$.% \cite[Chapter 4]{eisenbud20163264}.
\end{proof}

If we have a union of linear spaces $\Vcal = V_1 \cup \, \ldots, \, \cup V_\ell,$ the avoidance locus $\Acal_k(\Vcal)$ is the intersection of the avoidance loci $\Acal_k(V_1), \, \ldots, \, \Acal_k(V_l)$ by Lemma \ref{lem:easy_properties}. Thus the trichotomy in Lemma \ref{lem:linear} still holds.

For case (iii), we may compute the defining equations of these Schubert hyperplanes as follows. We parameterize the point $V$ as the kernel of a $(n-d) \times n$ matrix $Q,$ and a point $P$ in $\Gr(n-d,n)_\RR$ as the rowspan of a matrix $P.$ Then the linear space spanned by $P$ meets $V$ if and only if $\det(V \cdot P^T)$ vanishes. The Cauchy-Binet theorem then gives us a linear expression in the Pl\"ucker coordinates of $P.$

\begin{example}
Suppose $d = 0$ and $k = n-1$. We consider a collection of points $\Vcal = \{p_1,\ldots,p_\ell \} \subseteq \PP^{n-1}.$ Hyperplanes avoiding the points are points in the complements of hyperplanes $\{p_i\}^\perp$ in $(\PP^{n-1})^\vee$. The number of connected components in the avoidance locus $\Acal_{n-1}(\Vcal)$ is the number of regions in $(\PP^{n-1})^\vee \setminus \bigcup_{i=1}^{\ell}\{p_i\}^\perp$. 

The hyperplane arrangement $\bigcup_{i=1}^{\ell}\{p_i\}^\perp$ viewed in anl affine chart has an associated characteristic polynomial $\chi(t)$ \cite{stanley2007introduction}.
The number of projective regions in $(\PP^{n-1})^\vee \setminus \bigcup_{i=1}^{\ell}\{p_i\}^\perp$ is $\frac{\chi(1) + \chi(-1)}{2},$ where $\chi(1)$ is the number of bounded regions and $\chi(-1)$ is the number of all regions. For generic choices of $p_1,\ldots,p_\ell$ and affine chart, we have 
\begin{equation}\label{eq:regioncount}
   \chi(t) = t^{n-1} - \ell t^{n-2} + \binom{\ell} {2} t^{n-3} -\cdots + (-1)^{n-2}\binom{\ell}{n-1}. 
\end{equation}
\end{example}

It is difficult in general to compute the number of real regions of the complement of a Schubert arrangement in $\Gr(k,n)_\RR$. In partiular, this number is not determined by the characteristic polynomial of the hyperplanes in the ambient Pl\"ucker space. 

\begin{remark}
When the ambient Grassmannian is not a projective space, there is no general formula for the number of regions in the complement of hyperplanes.
For example, for six general Schubert hyperplanes in $\Gr(2,4)$, the number of regions may vary from 32 to 41; see \cite[Table 6]{mazzucchelli2025hyperplane}.
Moreover, a single region cut out by Schubert hyperplanes in the ambient Pl\"ucker space may intersect the Grassmannian in more than one region. Finally, regions may not be contractible in the ambient space, hence not convex; see \cite[Example 5.1,5.2]{mazzucchelli2025hyperplane}.
\end{remark}

\subsection{Curves}
The avoidance locus of an even-degree curve $C$ in $\bP^{n-1}$ consists of all the hyperplanes avoiding it. This is a union of regions in the dual variety of $C.$

Each region in $\cA_{n-1}(C)$ gives rise to a connected subset of $\cA_{n-2}(C)$ by Proposition \ref{prop:singleregion}. If we further assume that our curve is not a linear space, then $\text{conv}(C_i) \setminus C_i$ is nonempty for each component $C_i$ of $C$. Thus there are many points in $\cA_{n-2}(C)$ that do not arise via the incidence correspondence with $\cA_{n-1}(C)$. In particular, there is a Euclidean open set in $\Gr(n-2,n)$ of codimension two linear spaces passing through $\text{conv}(C_i) \setminus C_i,$ and a general one will not intersect $C_i.$

\begin{example}[A degree $6$ space curve]
Consider the degree six curve cut out by $x^2 - y^2 - xz = z -4x^3 + 3x = 0.$ This example appears in \cite{spacecurves}, in which they show it is a slice of a spectrahedron and compute its convex hull.
The avoidance locus $\mathcal{A}_1(X)$ is the complement of the Chow locus in $\Gr(2,4)_\RR.$ The Chow form is a polynomial of degree six with $35$ terms. Its complement has two components, since the real part of curve is homeomorphic to a circle. 

The Hurwitz form of the space curve is a homogeneous degree 10 equation with 62 terms. Its complement in $(\PP^3)^\vee$ has six regions, of which one region is the whole avoidance locus and the other five regions correspond to $2,4,4,4,$ and $6$ real intersection points. The avoidance locus yields a proper subset of one of the two regions in the complement of the dual variety.
\end{example}

\subsection{Hypersurfaces}

We have the tools to compute regions in the complement of a hypersurface in the affine space \cite{cummings2024smooth,breiding2025computing}. We present an approach to compute regions in the hypersurface complement of the Grassmannian $\Gr(2,n)_\RR$ in Algorithm \ref{alg:grassmannian}.

\begin{algorithm}[h]
\caption{Computing the regions in the complement of a hypersurface in~$\Gr(2,n)_\RR$} \label{alg:grassmannian}
\KwIn{Polynomial $f$ in Pl\"ucker coordinates $p_{01},\ldots,p_{n-2n-1}$ for the hypersurface.}

\KwOut{Regions in $\Gr(2,n)_\RR\setminus V(f)$, each with a sample point.}

Substitute the Pl\"ucker coordinates in $f$ by minors of $\begin{bmatrix}
1 & 0 & x_1 & \cdots & x_{n-2}\\
0 & 1 & y_1 & \cdots & y_{n-2}
\end{bmatrix}.$
Compute the affine regions in $\RR^{2n-4} \cong \Gr(2,n)_\RR\setminus V( p_{01})$.

Substitute the Pl\"ucker coordinates in $f$ by minors of $\begin{bmatrix}
1 & x_1 & 0 & x_2 & \cdots & x_{n-2}\\
0 & 0 & 1  & y_2 & \cdots & y_{n-2}
\end{bmatrix}.$
Compute the affine regions in $\RR^{2n-5}\cong V(p_{0,1})\setminus V(p_{0,2})$.

For each region result in step 2, we take a sample point $(x_1',x_2',\cdots, x_{n-2}', y_2',\cdots, y_{n-2}')$ and substitute the Pl\"ucker coordinates in $f$ by minors of 
$\begin{bmatrix}
1 & x_1' & 0 & x_2' & \cdots & x_{n-2}'\\
0 & t & 1  & y_2' & \cdots & y_{n-2}'
\end{bmatrix}$ to obtain $g(t)$. 

Solve $g(t)$ and pick $t_s>0$ such that $|t_s|$ is smaller than the absolute value of any real solution of $g(t) = 0$.

We form two points $p_ 1 = \begin{bmatrix}
1 & x_1' & 0 & x_2' & \cdots & x_{n-2}'\\
0 & t_s & 1  & y_2' & \cdots & y_{n-2}'
\end{bmatrix}, p_2 = \begin{bmatrix}
1 & x_1' & 0 & x_2' & \cdots & x_{n-2}'\\
0 & -t_s & 1  & y_2' & \cdots & y_{n-2}'
\end{bmatrix}$ and do a change of row basis to obtain $p_1',p_2'$ with the leftmost $2\times 2$ matrix the identity. 
The choice of $p_1',p_2'$ guarantee that they live in the same region of $\Gr(2,n)_\RR\setminus V(f)$ but are on two sides of $p_{01}=0$.

We track $p_1',p_2'$ to the regions they belong in the computation of Step 1. 
These regions are connected in $\Gr(2,n)_\RR\setminus V(f)$.
\end{algorithm}

\begin{example}[Quadric surface in $\PP^3$]
The topology of the real points of a quadric surface is controlled by the signature of the quadratic form. Signature $(4,0)$ means the variety is empty and $(3,1)$ means it may be written as a sphere in an affine chart. In the $(3,1)$ case, the avoidance loci $\cA_3(X)$ and $\cA_2(X)$ each consist of a single component, namely planes and lines outside the sphere.

Let us consider signature $(2,2)$; for instance, the surface given by $xz-yw=0.$ The complement of the Hurwitz form is a union of three connected components in $\Gr(2,4)_\RR$. The Hurwitz form has even degree, so it makes sense to talk about its sign. Two regions with $\Hcal_X<0$ form the avoidance locus $\Acal_2(X)$ while the third region has $\Hcal_X>0$. Lines in the third region intersects the surface in two points.
\end{example}

\begin{example}[Torus]
In Example \ref{eg:torushurwitz} we computed the Hurwitz form of the surface $X$ given by $(x^2+ y^2 + z^2 + 3 w^2)^2- 16(x^2+y^2)w^2 = 0.$ The complement of the Hurwitz locus in $\Gr(2,4)_{\RR}$ has four connected components. The avoidance locus consists of two of them; the other two consist of lines meeting $X$ in two or four real points, respectively. Topologically, the two regions in the avoidance locus correspond to lines passing through the torus hole and lines outside of it, see Figure \ref{fig:torus}.

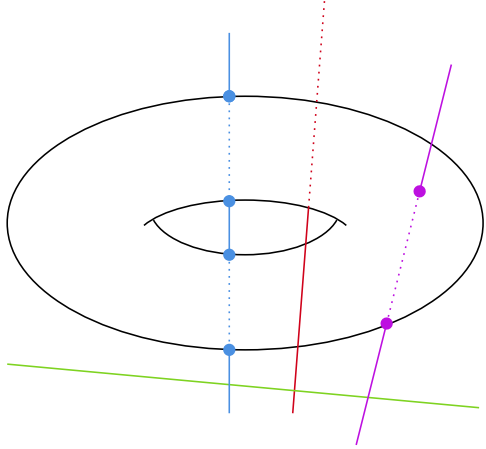
\begin{figure}[!h]
\begin{center}
    \scalebox{0.8}{\tikzset{every picture/.style={line width=0.75pt}} %set default line width to 0.75pt        

\begin{tikzpicture}[x=0.75pt,y=0.75pt,yscale=-1,xscale=1]

%Shape: Ellipse [id:dp38038586645836026] 
\draw   (160,140) .. controls (160,95.82) and (227.16,60) .. (310,60) .. controls (392.84,60) and (460,95.82) .. (460,140) .. controls (460,184.18) and (392.84,220) .. (310,220) .. controls (227.16,220) and (160,184.18) .. (160,140) -- cycle ;
%Shape: Arc [id:dp0028559740366069786] 
\draw  [draw opacity=0] (367.78,138.12) .. controls (360.69,150.74) and (337.51,160) .. (310,160) .. controls (282.49,160) and (259.31,150.74) .. (252.22,138.12) -- (310,130) -- cycle ; \draw   (367.78,138.12) .. controls (360.69,150.74) and (337.51,160) .. (310,160) .. controls (282.49,160) and (259.31,150.74) .. (252.22,138.12) ;  
%Shape: Arc [id:dp8963113828639865] 
\draw  [draw opacity=0] (246.22,141.5) .. controls (257.22,132.06) and (281.64,125.5) .. (310,125.5) .. controls (338.36,125.5) and (362.78,132.06) .. (373.78,141.5) -- (310,152.75) -- cycle ; \draw   (246.22,141.5) .. controls (257.22,132.06) and (281.64,125.5) .. (310,125.5) .. controls (338.36,125.5) and (362.78,132.06) .. (373.78,141.5) ;  
%Straight Lines [id:da8544300254996685] 
\draw [color={rgb, 255:red, 208; green, 2; blue, 27 }  ,draw opacity=1 ]   (350,130) -- (340,260) ;
%Straight Lines [id:da790267017419323] 
\draw [color={rgb, 255:red, 208; green, 2; blue, 27 }  ,draw opacity=1 ] [dash pattern={on 0.84pt off 2.51pt}]  (360,0) -- (350,130) ;
%Straight Lines [id:da8141763940899909] 
\draw [color={rgb, 255:red, 126; green, 211; blue, 33 }  ,draw opacity=1 ][fill={rgb, 255:red, 126; green, 211; blue, 33 }  ,fill opacity=1 ]   (457.5,256.5) -- (160,229) ;
%Straight Lines [id:da09157156284422041] 
\draw [color={rgb, 255:red, 189; green, 16; blue, 224 }  ,draw opacity=1 ]   (440,40) -- (420,120) ;
\draw [shift={(420,120)}, rotate = 104.04] [color={rgb, 255:red, 189; green, 16; blue, 224 }  ,draw opacity=1 ][fill={rgb, 255:red, 189; green, 16; blue, 224 }  ,fill opacity=1 ][line width=0.75]      (0, 0) circle [x radius= 3.35, y radius= 3.35]   ;
%Straight Lines [id:da4292480820080554] 
\draw [color={rgb, 255:red, 189; green, 16; blue, 224 }  ,draw opacity=1 ] [dash pattern={on 0.84pt off 2.51pt}]  (420,120) -- (399.2,203.4) ;
%Straight Lines [id:da42660534370894887] 
\draw [color={rgb, 255:red, 189; green, 16; blue, 224 }  ,draw opacity=1 ]   (399.2,203.4) -- (380,280) ;
\draw [shift={(399.2,203.4)}, rotate = 104.07] [color={rgb, 255:red, 189; green, 16; blue, 224 }  ,draw opacity=1 ][fill={rgb, 255:red, 189; green, 16; blue, 224 }  ,fill opacity=1 ][line width=0.75]      (0, 0) circle [x radius= 3.35, y radius= 3.35]   ;
%Straight Lines [id:da5773660053778377] 
\draw [color={rgb, 255:red, 74; green, 144; blue, 226 }  ,draw opacity=1 ]   (300,20) -- (300,60) ;
\draw [shift={(300,60)}, rotate = 90] [color={rgb, 255:red, 74; green, 144; blue, 226 }  ,draw opacity=1 ][fill={rgb, 255:red, 74; green, 144; blue, 226 }  ,fill opacity=1 ][line width=0.75]      (0, 0) circle [x radius= 3.35, y radius= 3.35]   ;
%Straight Lines [id:da7663437075367802] 
\draw [color={rgb, 255:red, 74; green, 144; blue, 226 }  ,draw opacity=1 ]   (300,220) -- (300,260) ;
%Straight Lines [id:da6557828222620147] 
\draw [color={rgb, 255:red, 74; green, 144; blue, 226 }  ,draw opacity=1 ] [dash pattern={on 0.84pt off 2.51pt}]  (300,60) -- (300,126.2) ;
\draw [shift={(300,126.2)}, rotate = 90] [color={rgb, 255:red, 74; green, 144; blue, 226 }  ,draw opacity=1 ][fill={rgb, 255:red, 74; green, 144; blue, 226 }  ,fill opacity=1 ][line width=0.75]      (0, 0) circle [x radius= 3.35, y radius= 3.35]   ;
\draw [shift={(300,60)}, rotate = 90] [color={rgb, 255:red, 74; green, 144; blue, 226 }  ,draw opacity=1 ][fill={rgb, 255:red, 74; green, 144; blue, 226 }  ,fill opacity=1 ][line width=0.75]      (0, 0) circle [x radius= 3.35, y radius= 3.35]   ;
%Straight Lines [id:da8200877510285283] 
\draw [color={rgb, 255:red, 74; green, 144; blue, 226 }  ,draw opacity=1 ] [dash pattern={on 0.84pt off 2.51pt}]  (300,160) -- (300,220) ;
\draw [shift={(300,220)}, rotate = 90] [color={rgb, 255:red, 74; green, 144; blue, 226 }  ,draw opacity=1 ][fill={rgb, 255:red, 74; green, 144; blue, 226 }  ,fill opacity=1 ][line width=0.75]      (0, 0) circle [x radius= 3.35, y radius= 3.35]   ;
\draw [shift={(300,160)}, rotate = 90] [color={rgb, 255:red, 74; green, 144; blue, 226 }  ,draw opacity=1 ][fill={rgb, 255:red, 74; green, 144; blue, 226 }  ,fill opacity=1 ][line width=0.75]      (0, 0) circle [x radius= 3.35, y radius= 3.35]   ;
%Straight Lines [id:da7591079095143134] 
\draw [color={rgb, 255:red, 74; green, 144; blue, 226 }  ,draw opacity=1 ]   (300,126.2) -- (300,160) ;

\end{tikzpicture}}
    \caption{A line in each of the four regions in the Hurwitz locus complement}\label{fig:torus}
\end{center}
\end{figure}

The avoidance locus $\Acal_3(X)$ of planes yields a subset of the second region by the incidence correspondence described in Proposition~\ref{prop:singleregion}. 
This subset is actually the whole second region.  
Indeed, consider a line $L$ tangent to $X$ and not passing through the interior of $\text{conv}(X)$. Since $X$ is compact in an affine chart of $\bP^3$, $\text{conv}(X)$ lies entirely on one side of any tangent plane. Thus we may choose a plane $H$ containing $L$ and perturb it slightly in the direction of the normal vector to obtain a plane in the avoidance locus of planes. 
\end{example}

\begin{example}[Quartic surface in $\PP^3$]
The Hurwitz form of the quartic surface $V(x_0^4 + x_1^4 - x_2^4 - x_3^4)$ has degree 12 and 68 monomials. The complement of the Hurwitz hypersurface in $\Gr(2,4)_\RR$ has 12 regions. Among them three regions form the avoidance locus, one region consists of linear spaces which intersect the surface in two real points, and the rest of the regions have linear spaces intersecting the surface in four real points.
\end{example}

\section{Positivity and Connected Components} \label{sec:n-1}
In this section we focus on the case $k = n-1,$ in which the avoidance locus lives in the dual projective space $(\PP^{n-1})^\vee.$ 
We give an upper bound for the number of regions in the avoidance locus in terms of the number of components of the real variety. 

The avoidance locus $\Acal_{n-1}(X)$ recovers certain known objects as special cases. In the introduction, we saw that for $X= v_d(\PP^{r-1})$ a Veronese embedding, the avoidance locus $\Acal_{n-1}(X)$ equals the open cone of positive polynomials with degree $d$ in $r$ variables. There is a generalization of this where we consider sparse positive polynomials with a specific support. 

\begin{definition}
    Let $A$ be a $d \times n$ matrix with integer entries, such that the vector $(1, 1, \, \ldots, \, 1)$ is in the rowspan of $A.$ The \emph{toric variety} $X_A$ is the closure of the image of
    \begin{align*}
        (\CC^*)^d & \to \PP^n \\
        t & \mapsto [t^{a_1}: \, \cdots \, :t^{a_r}].
    \end{align*}
\end{definition}

By the same argument as in the introduction, a hyperplane is in the avoidance locus $\Acal_{n-1}(X_A)$ if and only if the corresponding polynomial $\sum_{i = 1}^n c_i x^{a_i}$ is positive. 
This is called the \emph{sparse nonnegativity cone}, and approximations to it have been studied in \cite{forsgaard2022algebraic} as ways to certify the nonnegativity of sparse polynomials. 
This generalizes the theory of sums of squares, which certify nonnegativity of general polynomials, to \emph{sums of nonnegative circuits}, which certify nonnegativity of sparse polynomials.

We now prove an upper bound of regions in $\Acal_{n-1}(X)$.
The following argument is adapted from \cite[Proposition 5.3]{kaihnsa2019sixty} for the case of curves.

\begin{lemma}[Number of regions in $\Acal_{n-1}(X)$] \label{lem:regioncount}
    Let \( X \subset \bP^{n-1}\) be a smooth non-degenerate totally real variety of dimension $d$. Suppose that $X_\RR$ has $m$ components. Then, the number of regions in the avoidance locus is at most
    \[\frac{1}{4}\left(1 + m + \binom{m}{2} + \ldots + \binom{m}{n} + \binom{m-1}{n}\right).\]
\end{lemma}
\begin{proof}
    Let $X_1, \, \ldots, \, X_m$ be the connected components of $X_\RR.$ If the avoidance locus $\cA_{n-1}(X)$ is empty then we are done. Assume it is nonempty. Each hyperplane in the avoidance locus partitions these components into two sets. Since the variety is smooth and irreducible, each $X_i$ has the same dimension and no two components will meet; again, see \cite[Proposition 7.6.2]{bochnak2013real}. Now, take the convex hull of each $X_i$; if several convex hulls intersect, take the convex hull of the union of the components, and repeat this process until arriving at a collection of $r$ disjoint convex bodies, where $r \leq m$. Then each bipartition of these components corresponds uniquely to a region in the avoidance locus. The number of bipartitions of $r$ points in $\PP^{n-1}$ is equal to the number of regions in the complement of $r$ general hyperplanes. This correspondence is via projective duality between points and hyperplanes; the bipartitioning hyperplane becomes a point in $(\PP^{n-1})^\vee$, and the points on each side of it become hyperplanes in $(\PP^{n-1})^\vee$ whose defining linear forms are positive or negative on that point.
    
    The number of total regions and bounded regions of $r$ general hyperplanes in $\PP^{n-1}$ is given in Equation \eqref{eq:regioncount}, and is twice the upper bound in the theorem. The additional factor of $1/2$ comes from the non-adjacency of regions by Theorem \ref{thm:nonadjacent}.
\end{proof}

\begin{example}[Sparse nonnegativity cone]
    Any embedded toric variety $X_A$ has only one real connected component, because it has a real parameterization. This is consistent with the sparse nonnegativity cone being a single cone rather than a union of cones.
\end{example}

The upper bound in Lemma \ref{lem:regioncount} is quite coarse and can be refined by replacing $m$ with the number $r$ of components whose convex hulls are disjoint. The total number and topological type of the components is captured by the homology class of $X_\RR,$ and nesting behavior is captured by the isotopy type. A classification into isotopy type is known for plane curves of degree up to seven. A classification into isotopy types for surfaces in $\bP^3$ is known for degree up to four by work of Kharlamov \cite{kharlamov}.  For degree six plane curves, the known results and implications for the avoidance locus are explored in \cite[Section 5]{kaihnsa2019sixty}.

\section{Convexity of the avoidance locus}\label{sec:convexity}

In this section, we show several convexity properties of avoidance loci. We call a semialgebraic set $S \subset \bP^n$ \emph{convex} if the cone over $S$ is a convex cone in $\RR^{n+1},$ in the usual Euclidean sense. Our first result is that when the ambient Grassmannian is a projective space, each component of the avoidance locus is convex. 

\begin{theorem}\label{thm: hyperplane and X}
Let \( X \subset \bP^{n-1}\) be a smooth non-degenerate totally real variety of dimension $d$. The avoidance locus $\Acal_{n-1}(X) \subseteq (\PP^{n-1})^\vee$ is empty or a union of convex regions in $(\PP^{n-1}_\RR)^\vee \setminus L$ for some $L\in \Acal_{n-1}(X)$ up to closure.
\end{theorem}
\begin{proof}
We suppose the avoidance locus $\Acal_{n-1}(X)$ is not empty. Suppose $L\in\Acal_{n-1}(X)$, we will work in the affine chart $\PP^{n-1}\setminus L$.
We consider the connected components of $X_\RR$, denoted by $X_1,\ldots,X_k$ and their convex hull $\Conv(X_1),\ldots, \Conv(X_k)$. 
Any hyperplane intersecting $\Conv(X_i)$ intersects $X_i$ since if there is some $x = \sum_{j=1}^m \lambda_j x_j$ on the hyperplane with $\lambda_j \geq 0$, then either all $x_j$ lie on the hyperplane or there exist $x_j,x_j'$ that lie on two sides of the hyperplane. In the second case, we consider a continuous path from $x_j$ to $x_j'$ in $X_i$ and the path will pass through a point on $X_i$.
Hyperplanes not intersecting $\Conv(X_i)$ form the convex dual of $\Conv(X_i)$, denoted by $\Conv(X_i)^\vee$.
The avoidance locus in $\RR^n\cong (\PP^{n-1}_\RR)^\vee \setminus L$ is $\bigcap_{i=1}^k \Conv(X_i)$ taking away the points on the dual variety $X^\vee$, which is a union of disjoint convex sets up to closure.
\end{proof}

When the ambient Grassmannian is not a projective space, we must establish what we mean by convexity in a Grassmanian. 
There are multiple possible definitions, perhaps the simplest of which is convexity in the ambient Pl\"ucker space. This is known as \emph{extendable convexity} and was studied by Busemann  \cite{busemann}. However, in Example \ref{eg:notconvex} we demonstrate that the avoidance locus is not extendably convex in general.

\begin{example}\label{eg:notconvex}
Consider the union of two unit spheres with centers $(0,0,0)$ and $(2,2,0)$. The avoidance locus $\Acal_2(X)$ has one region since one can continuously move a line to avoid intersecting the two spheres. Slicing by the plane $E=\{z=0\}$ and restricting to $\Gr(2,E)\cong (\PP^2)^\vee$ produces the avoidance problem for two unit circles in $E$. 

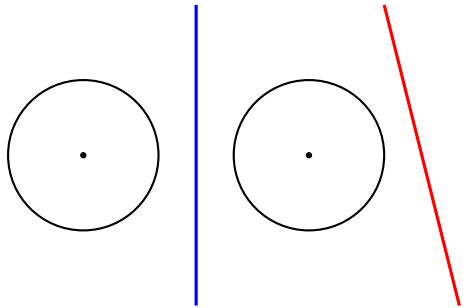
\begin{figure}[!h]
\centering
\begin{tikzpicture}[scale=1]
  % Two unit circles
  \draw[thick] (0,0) circle (1);
  \draw[thick] (3,0) circle (1);

  % Centers (optional markers)
  \fill (0,0) circle (1.2pt);
  \fill (3,0) circle (1.2pt);

  \draw[very thick, blue] (1.5,-2) -- (1.5,2) ;

  \draw[very thick, red] (5,-2) -- (4,2);
\end{tikzpicture}
\caption{Lines in the avoidance locus of two circles}\label{fig:circles}
\end{figure}
We computed that there are five regions in the complement of the dual curve in the plane, and two of them form the avoidance locus. Figure \ref{fig:circles} shows a line in each of these regions. One can see from the real picture or from Proposition \ref{prop:nonadjacent hyperplane} that these two regions are not adjacent, and hence the avoidance locus is not convex. In $\Gr(2,4),$ this translates to a line segment which connects the two lines from Figure \ref{fig:circles}, viewed as points in $\Gr(2,4).$ This segment is not entirely contained in the avoidance locus. 
Thus $\cA_2(X)$ is not convex in $\bP(\bigwedge^2\RR^4),$ and not extendably convex in $\Gr(2,4)_\RR.$
\end{example}

Instead, we use \emph{slice convexity}. This is a notion of convexity first defined by Shamovich and Vinnikov \cite{Shamovich_vinnikov} to study convexity of hyperbolicity regions in the Grassmanian. Let $E$ be a $k$-dimensional projective space. As in Definition \ref{def:schubert} we let $\Gr(k,E) \subset \Gr(k,n)$ denote the Schubert variety of all spaces of projective dimension $k-1$ contained in $E.$ In this context, we call $\Gr(k,E)$ a \emph{slice}.

\begin{definition}[Slice convex]
Let $E$ be a $k$-dimensional projective space. We say a connected set $S$ in $\Gr(k,n)_\RR$ is \emph{slice-convex with respect to $E$} if its intersection with the Schubert variety $\Gr(k,E)_\RR$ is convex in $\Gr(k,E)_\RR \cong \PP_\RR^k.$ 
We say that $S$ is \emph{slice-convex} if it is slice-convex with respect to all $E.$  We call $S$ \emph{component-wise slice-convex} with respect to $E$ if each connected component of $\Gr(k,E)_\RR \cap S$ is slice-convex.
\end{definition}

This generalizes the definition of convexity in projective space; indeed, a set $S$ is convex in projective space if its intersection with every line is convex in that line. We note that our setup slightly differs from that of \cite{Shamovich_vinnikov}. 
They say that $S$ is slice convex with respect to a $(k-2)$-dimensional linear space $E'$ if it is convex in the Schubert variety $\Gr(E', k)_\RR$ of linear spaces containing $E'.$  Taking the orthogonal complement everywhere recovers our definition.

\begin{theorem} \label{thm:slice convex}
Let \( X \subset \bP^{n-1}\) be a smooth non-degenerate totally real variety of dimension $d$. Then any avoidance locus of $X$ is component-wise slice-convex. 
\end{theorem}
\begin{proof}
Consider $\Acal_k(X)$ and let $E$ be a $k$-dimensional projective space. 
Let $\Acal_k(X\cap E,E)$ be all of the $(k-1)$-spaces contained in $E$ which intersect $X$ transversely in no real point.
The slice $\Acal_k(X)\cap \Gr(k,E)_\RR$ is equal to $\Acal_k(X\cap E,E)$ by Lemma \ref{lem:chowslice}. 
The result then follows from Theorem \ref{thm: hyperplane and X}. 
\end{proof}

\section*{Acknowledgements}
We are grateful to Bernd Sturmfels for helpful discussions and feedback throughout, and to Kristian Ranestad for useful comments on an early draft. Thanks also to Andrea Rosana, Rekha Thomas, and Kathlén Kohn for helpful discussions. We also thank the Max Planck Institute for Mathematics in the Sciences for providing the inspiring environment where this collaboration started. The first author was supported by NSF GRFP no. 2023358166.
\bibliographystyle{alpha}
\bibliography{reference}

@article{blekherman,
   title={Sums of squares and varieties of minimal degree},
   volume={29},
   ISSN={1088-6834},
   url={http://dx.doi.org/10.1090/jams/847},
   DOI={10.1090/jams/847},
   number={3},
   journal={Journal of the American Mathematical Society},
   publisher={American Mathematical Society (AMS)},
   author={Blekherman, Grigoriy and Smith, Gregory and Velasco, Mauricio},
   year={2015},
   month=sep, pages={893–913} }

@ARTICLE{holme,
  title     = "The geometric and numerical properties of duality in projective
               algebraic geometry",
  author    = "Holme, Audun",
  journal   = "Manuscripta Math.",
  publisher = "Springer Nature",
  volume    =  61,
  number    =  2,
  pages     = "145--162",
  year      =  1988,
  language  = "en"
}

@article{busemann,
  title={Convexity on {G}rassmann manifolds},
  author={Busemann, Herbert},
  journal={Enseign. Math.(2)},
  volume={7},
  year={1961}
}

@book{blekherman2012semidefinite,
  title={Semidefinite {O}ptimization and {C}onvex {A}lgebraic {G}eometry},
  author={Blekherman, Grigoriy and Parrilo, Pablo A and Thomas, Rekha R},
  year={2012},
  publisher={SIAM}
}

@article{gkzhyperdet,
title = {Hyperdeterminants},
journal = {Advances in Mathematics},
volume = {96},
number = {2},
pages = {226-263},
year = {1992},
issn = {0001-8708},
doi = {https://doi.org/10.1016/0001-8708(92)90056-Q},
url = {https://www.sciencedirect.com/science/article/pii/000187089290056Q},
author = {Gelfand, Israel M and Kapranov, Mikhail M and Zelevinsky, Andrei V}
}

@article{breidingnotices,
  author    = {Paul Breiding and Bernd Sturmfels and Sascha Timme},
  title     = {3264 Conics in a Second},
  journal   = {Notices of the American Mathematical Society},
  volume    = {67},
  number    = {1},
  pages     = {30--37},
  year      = {2020},
  doi       = {10.1090/noti2010},
}

@ARTICLE{kaji,
  title     = "On the duals of {S}egre varieties",
  author    = "Kaji, Hajime",
  journal   = "Geom. Dedicata",
  publisher = "Springer Science and Business Media LLC",
  volume    =  99,
  number    =  1,
  pages     = "221--229",
  year      =  2003,
  copyright = "https://www.springernature.com/gp/researchers/text-and-data-mining",
  language  = "en"
}

@article{holweck,
title = {Singularities of duals of {G}rassmannians},
author = {Frédéric Holweck},
journal = {Journal of Algebra},
volume = {337},
number = {1},
pages = {369-384},
year = {2011},
issn = {0021-8693},
doi = {https://doi.org/10.1016/j.jalgebra.2011.04.023},
url = {https://www.sciencedirect.com/science/article/pii/S0021869311002377}
}

@BOOK{Lee,
  title     = "Introduction to {S}mooth {M}anifolds",
  author    = "Lee, John M",
  publisher = "Springer",
  series    = "Graduate Texts in Mathematics",
  edition   =  2,
  year      =  2012,
  address   = "New York, NY",
  copyright = "https://www.springernature.com/gp/researchers/text-and-data-mining",
  language  = "en"
}

@ARTICLE{spacecurves,
  title     = "On the convex hull of a space curve",
  author    = "Ranestad, Kristian and Sturmfels, Bernd",
  journal   = "Adv. Geom.",
  publisher = "Walter de Gruyter GmbH",
  volume    =  12,
  number    =  1,
  pages     = "157--178",
  year      =  2012,
  language  = "en"
}

@ARTICLE{kharlamov,
  title     = "Isotopic types of nonsingular surfaces of fourth degree in {RP3}",
  author    = "Kharlamov, Viatcheslav Mikhailovich",
  journal   = "Funct. Anal. Its Appl.",
  publisher = "Springer Science and Business Media LLC",
  volume    =  12,
  number    =  1,
  pages     = "68--69",
  year      =  1978,
  language  = "en"
}

@article{Shamovich_vinnikov,
  title={Livsic-type determinantal representations and hyperbolicity},
  author={Shamovich, Eli and Vinnikov, Victor},
  journal={Advances in Mathematics},
  volume={329},
  pages={487--522},
  year={2018},
  publisher={Elsevier}
}

@book{bochnak2013real,
  title={Real {A}lgebraic {G}eometry},
  author={Bochnak, Jacek and Coste, Michel and Roy, Marie-Fran{\c{c}}oise},
  volume={36},
  year={2013},
  publisher={Springer Science \& Business Media}
}

@article{sturmfels2017hurwitz,
  title={The {H}urwitz form of a projective variety},
  author={Sturmfels, Bernd},
  journal={Journal of Symbolic Computation},
  volume={79},
  pages={186--196},
  year={2017},
  publisher={Elsevier}
}

@article{Kohn,
title = {Coisotropic hypersurfaces in {G}rassmannians},
journal = {Journal of Symbolic Computation},
volume = {103},
pages = {157-177},
year = {2021},
issn = {0747-7171},
doi = {https://doi.org/10.1016/j.jsc.2019.12.002},
author = {Kathlén Kohn}}

@article{
kaihnsa2019sixty,
  title={Sixty-four curves of degree six},
  author={Kaihnsa, Nidhi and Kummer, Mario and Plaumann, Daniel and Namin, Mahsa Sayyary and Sturmfels, Bernd},
  journal={Experimental Mathematics},
  volume={28},
  number={2},
  pages={132--150},
  year={2019},
  publisher={Taylor \& Francis}
}

@article{ranestad2024real,
  title={A Real Generalized Trisecant Trichotomy},
  author={Ranestad, Kristian and Seigal, Anna and Wang, Kexin},
  journal={arXiv preprint arXiv:2409.01356},
  year={2024}
}

@article{mazzucchelli2025hyperplane,
  title={Hyperplane arrangements in the grassmannian},
  author={Mazzucchelli, Elia and Pavlov, Dmitrii and Wang, Kexin},
  journal={Le Matematiche},
  volume={80},
  number={1},
  pages={387--408},
  year={2025}
}

@article{stanley2007introduction,
  title={An {I}ntroduction to {H}yperplane {A}rrangements},
  author={Stanley, Richard P and others},
  journal={Geometric {C}ombinatorics},
  volume={13},
  pages={389--496},
  year={2007}
}

@book{harris2013algebraic,
  title={Algebraic {G}eometry: {A} {F}irst {C}ourse},
  author={Harris, Joe},
  volume={133},
  year={2013},
  publisher={Springer Science \& Business Media}
}

@book{GKZ,
  author={Gelfand, Israel M and Kapranov, Mikhail M and Zelevinsky, Andrei V},
  title={Discriminants, Resultants, and Multidimensional Determinants},
  pages={13--47},
  year={1994},
  publisher={Springer}
}

@article{cummings2024smooth,
  title={Smooth connectivity in real algebraic varieties},
  author={Cummings, Joseph and Hauenstein, Jonathan D and Hong, Hoon and Smyth, Clifford D},
  journal={Numerical Algorithms},
  pages={1--22},
  year={2024},
  publisher={Springer}
}

@article{breiding2025computing,
  title={Computing arrangements of hypersurfaces},
  author={Breiding, Paul and Sturmfels, Bernd and Wang, Kexin},
  journal={Journal of Software for Algebra and Geometry},
  volume={15},
  number={1},
  pages={11--27},
  year={2025},
  publisher={Mathematical Sciences Publishers}
}

@article{forsgaard2022algebraic,
  title={The algebraic boundary of the sonc-cone},
  author={Forsg{\aa}rd, Jens and de Wolff, Timo},
  journal={SIAM Journal on Applied Algebra and Geometry},
  volume={6},
  number={3},
  pages={468--502},
  year={2022},
  publisher={SIAM}
}
\end{document}